\documentclass[review,onefignum,onetabnum]{siamonline250211}

\usepackage{lipsum}
\usepackage{amsfonts}
\usepackage{graphicx}
\usepackage{epstopdf}
\usepackage{algorithmic}
\ifpdf
  \DeclareGraphicsExtensions{.eps,.pdf,.png,.jpg}
\else
  \DeclareGraphicsExtensions{.eps}
\fi

\usepackage{enumitem}
\setlist[enumerate]{leftmargin=.5in}
\setlist[itemize]{leftmargin=.5in}

\newsiamremark{remark}{Remark}
\newsiamremark{hypothesis}{Hypothesis}
\crefname{hypothesis}{Hypothesis}{Hypotheses}
\newsiamthm{claim}{Claim}
\newsiamremark{fact}{Fact}
\crefname{fact}{Fact}{Facts}
\newsiamremark{assumption}{Assumption}
\crefname{assumption}{Assumption}{Assumption}
\nolinenumbers

\headers{Shadow price}{T. Rogala and Ł. Stettner}

\title{Discrete time shadow price revisited\thanks{Submitted to the editors DATE.
\funding{The work of the second author was supported by  NCN grant 2024/53/B/ST1/00703.}}}

\author{Tomasz Rogala\thanks{Cardinal Stefan Wyszynski University, Warsaw, Poland
  (\email{rogalatp@gmail.com}).}
\and Łukasz Stettner \thanks{Institute of Mathematics Polish Acad. Sci., Warsaw, Poland
  (\email{stettner@impan.pl}). }}

\usepackage{amsopn}

\DeclareMathOperator*{\esssup}{ess\,sup}

\ifpdf
\hypersetup{
  pdftitle={Shadow price revisited},
  pdfauthor={T. Rogala and Ł. Stettner}}
\fi

\DeclareMathOperator\closure{cl}
\DeclareMathOperator\sgn{sgn}

\begin{document}

\maketitle

\begin{abstract}
In the paper discrete time shadow price is constructed for the market with several assets with given bid and ask prices. Shadow price is the price such that the problem of optimal utility from terminal wealth on the market without transaction costs gives the same value function as in the case of bid and ask prices. In the paper we solve first static problem for two assets and then we construct the shadow price for dynamic model. Finally we present construction of the shadow price for several assets.
\end{abstract}



\section{Introduction}

Assume we are given a financial market consisting of $d$-assets defined on a filtered  probability space $(\Omega,{\cal F},({\cal F}_t),\mathbb{P})$  with bid $(\underline{S}_1(t),\ldots,\underline{S}_d(t))$ and ask $(\overline{S}_1(t),\ldots,\overline{S}_d(t))$ prices at discrete time $t$ such that $\underline{S}_i(t)\leq \overline{S}_i(t)$ and $\underline{S}_i(t)$, $\overline{S}_i(t)$ are ${\cal F}_t$ - measurable for $i\in \left\{1,\ldots,d\right\}$ and $t=0,1,\ldots,T$. Our financial position will be denoted by the vector $(x(t),y_1(t),\ldots,y_d(t))$ of amount on the bank account $x(t)$ and number of $i$th assets in our portfolio $y_i(t)$ at time $t$, with $i=1,2,\ldots, d$ and 
$t=0,1,\ldots,T$. In what follows we assume that short selling and short borrowing is not allowed so that vector $(x(t),y_1(t),\ldots,y_d(t))$ is nonnegative. Given time horizon $T$ and utility function $U$, which is increasing and strictly concave we want to maximize utility from terminal wealth $W_T$ i.e.
\begin{equation}
\mathbb{E}\left\{U(W_T)\right\}
\end{equation}
where $W_T=x(T)+\sum_{i=1}^d y_i(T)\underline{S}_i(T)$.

The problem can be solved by backward induction following ideas of the papers \cite{RS1} and \cite{RS2}.
The purpose of this paper is to show existence of so called shadow price. It is a price on the market without transaction costs (bid and ask prices are the same) for which optimal value of utility from terminal wealth is the same as in the market with bid and ask price. It is not obvious that such price exists (see a counterexample in \cite{CziScha}). In mathematics of finance literature such problems appeared first in the paper \cite{KaMuh1} for the problems on the market with transaction costs and in \cite{Kuhn} for an investor trading in a limit order market.
Discrete time market model was studied in the papers \cite{KaMuh2},  \cite{Rok}, \cite{CziScha}, \cite{RS1}, \cite{RS2}. The case with finite probability space was considered in \cite{KaMuh2}. In \cite{Rok} game theory interpretation of shadow price was shown. In \cite{CziScha} shadow price was studied using duality approach. The papers  \cite{RS1}, \cite{RS2} construct shadow price using direct backward induction approach. In the paper  \cite{RS1} the case of one asset was considered. The paper \cite{RS2} is generalizing  \cite{RS1} to several-asset however complete results are only shown in two-asset case.  In this paper we introduce a methodology which allows us to consider several-asset case. Although to simplify notations we formulate our results for two asset case our methodology can be easily extended to multidimensional case as is shown in section 10.   The key argument introduced in this paper is based on the detailed analysis of the static case. Some properties of buying and selling zones are of independent interest. The approach to shadow price formulated respectively in this paper and \cite{RS1}, \cite{RS2} seems to be most natural and does not require strong assumptions. 

Continuous time shadow price was introduced first for lognormal asset prices in \cite{KaMuh1} and then continued in \cite{Ben}, \cite{Choi}, \cite{Ger}, \cite{Czi1}, \cite{Czi2}, \cite{Czi}, \cite{Czi3}, \cite{Dol}. The construction of shadow price was first studied for geometric Brownian motion asset prices and then extended to semimartingales and more general processes including geometric fractional Brownian motion prices, continuous processes or even c\`adl\`ag prices. The main tool is based on duality arguments.

The paper is organized as follows. We consider first shadow price for the static case with $t=1$ and two assets. Then we construct shadow price in the dynamic two asset case. The case of several assets can be considered by analogy to the two assets case and we point out only main steps to avoid rather cumbersome notations.
What is interesting in the paper is that probabilistic problem is in major part  reduced to studies of  geometry of buying, selling and no transaction zones in the static case. This part is extended beyond the scope necessary to construct shadow price, since we were able to show some interesting properties of the static buying and selling zones which in our opinion are of independent interest. Furthermore construction of multidimensional shadow price is based on the construction of shadow prices for each asset separately. Dynamical case studied in section 9 requires nondegeneracy of asset prices and is shown using additional nontrivial reasoning. 

\section{Static problem - preliminaries}

In this part of the paper we consider stationary two-asset case with bid $(\underline{s}_1,\underline{s}_2)$ and ask $(\overline{s}_1,\overline{s}_2)$ prices. We denote by $\mathbb{D}$ the set of prices
\begin{equation}
    \mathbb{D} := \big\{ (\underline{s}_{1}, \underline{s}_{2}, \overline{s}_{1}, \overline{s}_{2}) \in \mathbb{R}_{+}^{2} \times \mathbb{R}_{+}^{2} : \quad 0 < \underline{s}_{1} \leq \overline{s}_{1} \quad \mbox{and} \quad 0 < \underline{s}_{2} \leq \overline{s}_{2} \big\},
\end{equation}
where $\mathbb{R}_{+}=[0,\infty)$.
Our financial position will be denoted by $(x, y_{1}, y_{2})$ and will be nonnegative.
For every $(x, y_{1}, y_{2}, \underline{s}_{1}, \underline{s}_{2}, \overline{s}_{1}, \overline{s}_{2}) \in \mathbb{R}_{+} \times \mathbb{R}_{+}^{2} \times \mathbb{D}$ we have the following set of admissible investment strategies
\begin{eqnarray}
&&\mathbb{A}(x, y_{1}, y_{2}, \underline{s}_{1}, \underline{s}_{2}, \overline{s}_{1}, \overline{s}_{2}) := \nonumber \\
&& \Big\{ (l_{1}, l_{2}, m_{1}, m_{2}) \in \mathbb{R}_{+}^{2} \times \mathbb{R}_{+}^{2} : \quad 0 \leq x + \underline{s}_{1} m_{1} + \underline{s}_{2} m_{2} - \overline{s}_{1} l_{1} - \overline{s}_{2} l_{2}, \nonumber \\
&& 0 \leq y_{1} - m_{1} + l_{1}, 0 \leq y_{2} - m_{2} + l_{2}, m_{1} \leq y_{1} \quad \mbox{and} \quad m_{2} \leq y_{2} \Big\} .
\end{eqnarray}
We clearly have that

\begin{lemma}\label{Lem1.1}
    For every $(x, y_{1}, y_{2}, \underline{s}_{1}, \underline{s}_{2}, \overline{s}_{1}, \overline{s}_{2}) \in \mathbb{R}_{+} \times \mathbb{R}_{+}^{2} \times \mathbb{D}$ the set $\mathbb{A}(x, y_{1}, y_{2}, \underline{s}_{1}, \underline{s}_{2}, \overline{s}_{1}, \overline{s}_{2})$ is compact and convex subset of $\mathbb{R}_{+}^{2} \times \mathbb{R}_{+}^{2}$.
\end{lemma}

Moreover following the proof of Theorem 2.1 of \cite{RS1} we have

\begin{proposition}\label{propH}
 The mapping
        \begin{equation}
            (x, y_{1}, y_{2}, \underline{s}_{1}, \underline{s}_{2}, \overline{s}_{1}, \overline{s}_{2}) \longmapsto \mathbb{A}(x, y_{1}, y_{2}, \underline{s}_{1}, \underline{s}_{2}, \overline{s}_{1}, \overline{s}_{2})
        \end{equation}
        is continuous in the Hausdorff metric.
\end{proposition}

Next three Lemmas show several simple properties of the sets $\mathbb{A}$

\begin{lemma}\label{Lem1.3}
    Let $(\underline{s}_{1}^{1}, \underline{s}_{2}^{1}, \overline{s}_{1}^{1}, \overline{s}_{2}^{1}), (\underline{s}_{1}^{2}, \underline{s}_{2}^{2}, \overline{s}_{1}^{2}, \overline{s}_{2}^{2}) \in \mathbb{D}$ be such that
    \begin{equation*}
        0 < \underline{s}_{1}^{1} \leq \underline{s}_{1}^{2} \leq \overline{s}_{1}^{2} \leq \overline{s}_{1}^{1} \quad \mbox{and} \quad 0 < \underline{s}_{2}^{1} \leq \underline{s}_{2}^{2} \leq \overline{s}_{2}^{2} \leq \overline{s}_{2}^{1} .
    \end{equation*}
    Then for every $(x, y_{1}, y_{2}) \in \mathbb{R}_{+} \times \mathbb{R}_{+}^{2}$ we have that
    \begin{equation}\label{inclusion for better situation}
        \mathbb{A} (x, y_{1}, y_{2}, \underline{s}_{1}^{1}, \underline{s}_{2}^{1}, \overline{s}_{1}^{1}, \overline{s}_{2}^{1}) \subseteq \mathbb{A} (x, y_{1}, y_{2}, \underline{s}_{1}^{2}, \underline{s}_{2}^{2}, \overline{s}_{1}^{2}, \overline{s}_{2}^{2}) .
    \end{equation}
\end{lemma}

\begin{lemma}\label{Lem1.4}
    Let $(x, y_{1}, y_{2}, \underline{s}_{1}, \underline{s}_{2}, \overline{s}_{1}, \overline{s}_{2}) \in \mathbb{R}_{+} \times \mathbb{R}_{+}^{2} \times \mathbb{D}$ and
    \begin{equation*}
        (l_{1}, l_{2}, m_{1}, m_{2}) \in \mathbb{A}(x, y_{1}, y_{2}, \underline{s}_{1}, \underline{s}_{2}, \overline{s}_{1}, \overline{s}_{2}).
    \end{equation*}
    Then when $l_{1} = m_{1}$ and $\widetilde{s} \in [\underline{s}_{1}, \overline{s}_{1}]$ we have
    \begin{equation}\label{wasting the first stock with the same price}
        (0, l_{2}, 0, m_{2}) \in \mathbb{A}(x, y_{1}, y_{2}, \widetilde{s}, \underline{s}_{2}, \widetilde{s}, \overline{s}_{2}),
    \end{equation}
    while when  $l_{2} = m_{2}$ and $\widetilde{s} \in [\underline{s}_{2}, \overline{s}_{2}]$ we have
    \begin{equation}\label{wasting the second stock with the same price}
        (l_{1}, 0, m_{1}, 0) \in \mathbb{A}(x, y_{1}, y_{2}, \underline{s}_{1}, \widetilde{s}, \overline{s}_{1}, \widetilde{s}).
    \end{equation}
\end{lemma}

\begin{lemma}\label{Lem1.5}
 Let $(x, y_{1}, y_{2}, \underline{s}_{1}, \underline{s}_{2}, \overline{s}_{1}, \overline{s}_{2}) \in \mathbb{R}_{+} \times \mathbb{R}_{+}^{2} \times \mathbb{D}$ and  $\widetilde{s} \in [\underline{s}_{1}, \overline{s}_{1}]$. Then when
    \begin{equation*}
        (l_{1}, l_{2}, m_{1}, m_{2}) \in \mathbb{A}(x, y_{1}, y_{2}, \widetilde{s}, \underline{s}_{2}, \widetilde{s}, \overline{s}_{2})
    \end{equation*}
    and $l_{1} = m_{1}$ we have
      \begin{equation}\label{admissibility when wasting on the first shadow stock}
        (0, l_{2}, 0, m_{2}) \in \mathbb{A}(x, y_{1}, y_{2}, \underline{s}_{1}, \underline{s}_{2}, \overline{s}_{1}, \overline{s}_{2}) .
    \end{equation}
In the case when
\begin{equation*}
        (l_{1}, l_{2}, m_{1}, m_{2}) \in \mathbb{A}(x, y_{1}, y_{2}, \underline{s}_{1}, \widetilde{s}, \overline{s}_{1}, \widetilde{s})
    \end{equation*}
     and $l_{2} = m_{2}$ we have
     \begin{equation}\label{admissibility when wasting on the second shadow stock}
        (l_{1}, 0, m_{1}, 0) \in \mathbb{A}(x, y_{1}, y_{2}, \underline{s}_{1}, \underline{s}_{2}, \overline{s}_{1}, \overline{s}_{2}).
\end{equation}
\end{lemma}

The following Proposition will be important in further constructions

\begin{proposition}\label{Prop1.6}
    Let $(x, y_{1}, y_{2}, \underline{s}_{1}, \underline{s}_{2}, \overline{s}_{1}, \overline{s}_{2}) \in \mathbb{R}_{+} \times \mathbb{R}_{+}^{2} \times \mathbb{D}$,
    \begin{equation*}
        (\widehat{l}_{1}, \widehat{l}_{2}, \widehat{m}_{1}, \widehat{m}_{2}) \in \mathbb{A}(x, y_{1}, y_{2}, \underline{s}_{1}, \underline{s}_{2}, \overline{s}_{1}, \overline{s}_{2})
    \end{equation*}
    and
    \begin{equation*}
        (\widetilde{l}_{1}, \widetilde{l}_{2}, \widetilde{m}_{1}, \widetilde{m}_{2}) \in \mathbb{A}(x, y_{1}, y_{2}, \overline{s}_{1}, \underline{s}_{2}, \overline{s}_{1}, \overline{s}_{2}).
    \end{equation*}
    Then when
    \begin{equation*}
        \widehat{l}_{1} > 0, \widehat{m}_{1} = 0 \quad \mbox{and} \quad \widetilde{l}_{1} = 0, \widetilde{m}_{1} > 0.
    \end{equation*}
we have
    \begin{equation}\label{admissibility of mixed strategies with upper shadow for the first stock}
        \Big( 0, \lambda \widehat{l}_{2} + (1 - \lambda ) \widetilde{l}_{2}, 0, \lambda \widehat{m}_{2} + (1 - \lambda ) \widetilde{m}_{2} \Big) \in \mathbb{A}(x, y_{1}, y_{2}, \underline{s}_{1}, \underline{s}_{2}, \overline{s}_{1}, \overline{s}_{2}) .
    \end{equation}
    with
    \begin{equation*}
        \lambda := \frac{\widetilde{m}_{1}}{\widehat{l}_{1} + \widetilde{m}_{1}}.
        \end{equation*}
        Similarly, for
        \begin{equation*}
        (\widehat{l}_{1}, \widehat{l}_{2}, \widehat{m}_{1}, \widehat{m}_{2}) \in \mathbb{A}(x, y_{1}, y_{2}, \underline{s}_{1}, \underline{s}_{2}, \overline{s}_{1}, \overline{s}_{2})
    \end{equation*}
    and
    \begin{equation*}
        (\widetilde{l}_{1}, \widetilde{l}_{2}, \widetilde{m}_{1}, \widetilde{m}_{2}) \in \mathbb{A}(x, y_{1}, y_{2}, \underline{s}_{1}, \overline{s}_{2}, \overline{s}_{1}, \overline{s}_{2})
    \end{equation*}
    such that
    \begin{equation*}
        \widehat{l}_{2} > 0, \widehat{m}_{2} = 0 \quad \mbox{and} \quad \widetilde{l}_{2} = 0, \widetilde{m}_{2} > 0.
    \end{equation*}
    we have
    \begin{equation}\label{admissibility of mixed strategies with upper shadow for the second stock}
        \Big( \lambda \widehat{l}_{1} + (1 - \lambda ) \widetilde{l}_{1}, 0, \lambda \widehat{m}_{1} + (1 - \lambda ) \widetilde{m}_{1}, 0 \Big) \in \mathbb{A}(x, y_{1}, y_{2}, \underline{s}_{1}, \underline{s}_{2}, \overline{s}_{1}, \overline{s}_{2}) .
    \end{equation}
     with \begin{equation*}
        \lambda := \frac{\widetilde{m}_{2}}{\widehat{l}_{2} + \widetilde{m}_{2}}.
        \end{equation*}
\end{proposition}
\begin{proof} We show only that for $\lambda=\frac{\widetilde{m}_{1}}{\widehat{l}_{1} + \widetilde{m}_{1}}$ we have
     $(\ref{admissibility of mixed strategies with upper shadow for the first stock})$.
The proof of the second part of Proposition is similar.
     Since $\underline{s}_{1} \leq \overline{s}_{1} \leq \overline{s}_{1} \leq \overline{s}_{1}$ and $(\widehat{l}_{1}, \widehat{l}_{2}, \widehat{m}_{1}, \widehat{m}_{2}) \in \mathbb{A}(x, y_{1}, y_{2}, \underline{s}_{1}, \underline{s}_{2}, \overline{s}_{1}, \overline{s}_{2})$ by  Lemma $\ref{Lem1.3}$ we have
        $(\widehat{l}_{1}, \widehat{l}_{2}, \widehat{m}_{1}, \widehat{m}_{2}) \in \mathbb{A}(x, y_{1}, y_{2}, \overline{s}_{1}, \underline{s}_{2}, \overline{s}_{1}, \overline{s}_{2})$.
     By Lemma $\ref{Lem1.1}$ the set $\mathbb{A}(x, y_{1}, y_{2}, \overline{s}_{1}, \underline{s}_{2}, \overline{s}_{1}, \overline{s}_{2})$ is convex and therefore
    \begin{equation*}
        \begin{split}
            &\lambda (\widehat{l}_{1}, \widehat{l}_{2}, \widehat{m}_{1}, \widehat{m}_{2}) + (1 - \lambda ) (\widetilde{l}_{1}, \widetilde{l}_{2}, \widetilde{m}_{1}, \widetilde{m}_{2}) = \\
            &=\Big( \lambda \widehat{l}_{1} + (1 - \lambda ) \widetilde{l}_{1}, \lambda \widehat{l}_{2} + (1 - \lambda ) \widetilde{l}_{2}, \lambda \widehat{m}_{1} + (1 - \lambda ) \widetilde{m}_{1}, \lambda \widehat{m}_{2} + (1 - \lambda ) \widetilde{m}_{2} \Big) = \\
            &=\Big( \lambda \widehat{l}_{1}, \lambda \widehat{l}_{2} + (1 - \lambda ) \widetilde{l}_{2}, (1 - \lambda ) \widetilde{m}_{1}, \lambda \widehat{m}_{2} + (1 - \lambda ) \widetilde{m}_{2} \Big) = \\
            &=\Bigg( \frac{\widehat{l}_{1} \widetilde{m}_{1}}{\widehat{l}_{1} + \widetilde{m}_{1}}, \lambda \widehat{l}_{2} + (1 - \lambda ) \widetilde{l}_{2}, \frac{\widehat{l}_{1} \widetilde{m}_{1}}{\widehat{l}_{1} + \widetilde{m}_{1}}, \lambda \widehat{m}_{2} + (1 - \lambda ) \widetilde{m}_{2} \Bigg)
        \end{split}
    \end{equation*}
    is an element of
    \begin{equation*}
        \Bigg( \frac{\widehat{l}_{1} \widetilde{m}_{1}}{\widehat{l}_{1} + \widetilde{m}_{1}}, \lambda \widehat{l}_{2} + (1 - \lambda ) \widetilde{l}_{2}, \frac{\widehat{l}_{1} \widetilde{m}_{1}}{\widehat{l}_{1} + \widetilde{m}_{1}}, \lambda \widehat{m}_{2} + (1 - \lambda ) \widetilde{m}_{2} \Bigg) \in \mathbb{A}(x, y_{1}, y_{2}, \overline{s}_{1}, \underline{s}_{2}, \overline{s}_{1}, \overline{s}_{2}).
    \end{equation*}
    Consequently by Lemma $\ref{Lem1.5}$ we have that  $(\ref{admissibility of mixed strategies with upper shadow for the first stock})$.
\end{proof}

Lemma below shows us important property of the set of admissible strategies $\mathbb{A}$ which allows us to eliminate strategies based on simultaneous buying and selling of the same assets.
\begin{lemma}\label{Lem1.7}
    Let $(x, y_{1}, y_{2}, \underline{s}_{1}, \underline{s}_{2}, \overline{s}_{1}, \overline{s}_{2}) \in \mathbb{R}_{+} \times \mathbb{R}_{+} \times \mathbb{D}$ and assume that $(l_{1}, l_{2}, m_{1}, m_{2}) \in \mathbb{A}(x, y_{1}, y_{2}, \underline{s}_{1}, \underline{s}_{2}, \overline{s}_{1}, \overline{s}_{2})$. Define
    \begin{equation}\label{changing a strategy to an efficient strategy}
        (\widehat{l}_{1}, \widehat{l}_{2}, \widehat{m}_{1}, \widehat{m}_{2}) := \Big( (l_{1} - m_{1})^{+}, (l_{2} - m_{2})^{+}, (m_{1} - l_{1})^{+}, (m_{2} - l_{2})^{+} \Big) .
    \end{equation}
    Then $(\widehat{l}_{1}, \widehat{l}_{2}, \widehat{m}_{1}, \widehat{m}_{2}) \in \mathbb{A} (x, y_{1}, y_{2}, \underline{s}_{1}, \underline{s}_{2}, \overline{s}_{1}, \overline{s}_{2})$ and $\widehat{l}_{1} \widehat{m}_{1} + \widehat{l}_{2} \widehat{m}_{2} = 0$. Moreover,
    \begin{equation*}
        x + \underline{s}_{1} m_{1} + \underline{s}_{2} m_{2} - \overline{s}_{1} l_{1} - \overline{s}_{2} l_{2} \leq x + \underline{s}_{1} \widehat{m}_{1} + \underline{s}_{2} \widehat{m}_{2} - \overline{s}_{1} \widehat{l}_{1} - \overline{s}_{2} \widehat{l}_{2}
    \end{equation*}
    and
    \begin{equation*}
        y_{1} - m_{1} + l_{1} = y_{1} - \widehat{m}_{1} + \widehat{l}_{1} \quad \mbox{and} \quad y_{2} - m_{2} + l_{2} = y_{2} - \widehat{m}_{2} + \widehat{l}_{2} .
    \end{equation*}
\end{lemma}

\begin{proof}
    From Lemma $\ref{Lem1.2}$ we have that
    \begin{equation*}
        \begin{split}
            0 &\leq x + \underline{s}_{1} m_{1} + \underline{s}_{2} m_{2} - \overline{s}_{1} l_{1} - \overline{s}_{2} l_{2} = x + (\underline{s}_{1} m_{1} - \overline{s}_{1} l_{1}) + (\underline{s}_{2} m_{2}  - \overline{s}_{2} l_{2}) \leq \\
            &\leq x + \big[ \underline{s}_{1} (m_{1} - l_{1})^{+} - \overline{s}_{1} (l_{1} - m_{1})^{+} \big] + \big[ \underline{s}_{1} (m_{1} - l_{1})^{+} - \overline{s}_{1} (l_{1} - m_{1})^{+} \big] = \\
            &= x + (\underline{s}_{1} \widehat{m}_{1} - \overline{s}_{1} \widehat{l}_{1} )+ (\underline{s}_{2} \widehat{m}_{2} - \overline{s}_{2} \widehat{l}_{2}) = x + \underline{s}_{1} \widehat{m}_{1} + \underline{s}_{2} \widehat{m}_{2} - \overline{s}_{1} \widehat{l}_{1} - \overline{s}_{2} \widehat{l}_{2} .
        \end{split}
    \end{equation*}
    Clearly, for $i = 1, 2$ we have that
    \begin{equation*}
        \begin{split}
            0 & \leq y_{i} - m_{i} + l_{i} = y_{i} - (m_{i} - l_{i}) = \big[ (m_{i} - l_{i})^{+} - (l_{i} - m_{i})^{+} \big] = \\
        & = y_{i} - (\widehat{m}_{i} - \widehat{l}_{i}) = y_{i} - \widehat{m}_{i} + \widehat{l}_{i} .
        \end{split}
    \end{equation*}
    and $\widehat{l}_{1} \widehat{m}_{1} + \widehat{l}_{2} \widehat{m}_{2} = 0$.
\end{proof}

\section{Static value function}
Let
\begin{equation}
    \begin{split}
        \mathbb{P} := \Big\{ (x, y_{1}, y_{2}, \underline{s}_{1}, \underline{s}_{2}, \overline{s}_{1}, \overline{s}_{2}, l_{1}, l_{2}, & m_{1}, m_{2}) \in \mathbb{R}_{+} \times \mathbb{R}_{+}^{2} \times \mathbb{D} \times \mathbb{R}_{+}^{2} \times \mathbb{R}_{+}^{2}:  \\
        & (l_{1}, l_{2}, m_{1}, m_{2}) \in \mathbb{A}(x, y_{1}, y_{2}, \underline{s}_{1}, \underline{s}_{2}, \overline{s}_{1}, \overline{s}_{2}) \Big\} .
    \end{split}
\end{equation}

Let $g : \mathbb{R}_{+} \times \mathbb{R}_{+}^{2} \longrightarrow \mathbb{R} \cup \{ - \infty \}$ be a strictly concave function which is strictly increasing with respect to each variable. Such function appears in a natural way when we consider optimal utility from terminal wealth with strictly concave utility function (see Lemma \ref{concaveutil}).

Define the function $\varphi : \mathbb{P} \longrightarrow \mathbb{R} \cup \{ - \infty \}$ by the formula
\begin{eqnarray}\label{vphi}
&& \varphi (x, y_{1}, y_{2}, \underline{s}_{1}, \underline{s}_{2}, \overline{s}_{1}, \overline{s}_{2}, l_{1}, l_{2}, m_{1}, m_{2}) := \nonumber \\
&& g \Big( x + \underline{s}_{1} m_{1} + \underline{s}_{2} m_{2} - \overline{s}_{1} l_{1} - \overline{s}_{2} l_{2}, y_{1} - m_{1} + l_{1}, y_{2} - m_{2} + l_{2} \Big)
\end{eqnarray}
and
$\Phi : \mathbb{R}_{+} \times \mathbb{R}_{+}^{2} \times \mathbb{D} \longrightarrow \mathbb{R} \cup \{ - \infty \}$ for every $(x, y_{1}, y_{2}, \underline{s}_{1}, \underline{s}_{2}, \overline{s}_{1}, \overline{s}_{2}) \in \mathbb{R}_{+} \times \mathbb{R}_{+}^{2} \times \mathbb{D}$
\begin{eqnarray}\label{Phi}
&&\Phi (x, y_{1}, y_{2}, \underline{s}_{1}, \underline{s}_{2}, \overline{s}_{1}, \overline{s}_{2}) := \nonumber \\
&&\sup_{(l_{1}, l_{2}, m_{1}, m_{2}) \in \mathbb{A}(x, y_{1}, y_{2}, \underline{s}_{1}, \underline{s}_{2}, \overline{s}_{1}, \overline{s}_{2})} \varphi (x, y_{1}, y_{2}, \underline{s}_{1}, \underline{s}_{2}, \overline{s}_{1}, \overline{s}_{2}, l_{1}, l_{2}, m_{1}, m_{2}) .
\end{eqnarray}

\begin{lemma}\label{Lem2.01}
Functions $\varphi$ and $\Phi$ are concave for fixed  $(\underline{s}_{1}, \underline{s}_{2}, \overline{s}_{1}, \overline{s}_{2})$  and continuous with respect to its variables. Moreover there is an element $(l_{1}, l_{2}, m_{1}, m_{2}) \in $ $\mathbb{A}(x, y_{1}, y_{2}, \underline{s}_{1}, \underline{s}_{2}, \overline{s}_{1}, \overline{s}_{2})$ for which supremum in  \eqref{Phi} is attained.

\begin{proof}
Since $\varphi$ is a composition of a concave function $g$ and affine functions it is concave for fixed  $(\underline{s}_{1}, \underline{s}_{2}, \overline{s}_{1}, \overline{s}_{2})$. Moreover $g$ as a concave function is continuous and $\varphi$ is also continuous with respect to its all variables. Furthermore since the set $\mathbb{A}(x, y_{1}, y_{2}, \underline{s}_{1}, \underline{s}_{2}, \overline{s}_{1}, \overline{s}_{2})$ is compact and continuous in Hausdorff metric with respect to $(x, y_{1}, y_{2},\underline{s}_{1}, \underline{s}_{2}, \overline{s}_{1}, \overline{s}_{2})$ (see Proposition \ref{propH}), we easily have continuity of $\Phi$ and existence of an element for which supremum is attained.
\end{proof}

\begin{lemma}\label{Lem2.1}
    Let $(x, y_{1}, y_{2}, \underline{s}_{1}, \underline{s}_{2}, \overline{s}_{1}, \overline{s}_{2}, l_{1}, l_{2}, m_{1}, m_{2}) \in \mathbb{P}$ and let $(\widehat{l}_{1}, \widehat{l}_{2}, \widehat{m}_{1}, \widehat{m}_{2})$ be defined by $(\ref{changing a strategy to an efficient strategy})$.
    Then
    \begin{equation*}
        \varphi (x, y_{1}, y_{2}, \underline{s}_{1}, \underline{s}_{2}, \overline{s}_{1}, \overline{s}_{2}, l_{1}, l_{2}, m_{1}, m_{2}) \leq \varphi (x, y_{1}, y_{2}, \underline{s}_{1}, \underline{s}_{2}, \overline{s}_{1}, \overline{s}_{2}, \widehat{l}_{1}, \widehat{l}_{2}, \widehat{m}_{1}, \widehat{m}_{2}) .
    \end{equation*}
\end{lemma}

\begin{proof}
By Lemma \ref{Lem1.7} clearly $(\widehat{l}_{1}, \widehat{l}_{2}, \widehat{m}_{1}, \widehat{m}_{2}) \in \mathbb{A} (x, y_{1}, y_{2}, \underline{s}_{1}, \underline{s}_{2}, \overline{s}_{1}, \overline{s}_{2})$.
Since function $g$ is strictly increasing with respect to each variable, by Lemma \ref{Lem1.7} again  we have
    \begin{eqnarray*}
&&\varphi (x, y_{1}, y_{2}, \underline{s}_{1}, \underline{s}_{2}, \overline{s}_{1}, \overline{s}_{2},
 l_{1}, l_{2}, m_{1}, m_{2}) =  \nonumber \\
 && g ( x + \underline{s}_{1} m_{1} + \underline{s}_{2} m_{2} - \overline{s}_{1} l_{1} - \overline{s}_{2} l_{2}, y_{1} - m_{1} + l_{1}, y_{2} - m_{1} + l_{2}) \leq  \nonumber \\
&& \leq  g ( x + \underline{s}_{1} \widehat{m}_{1} + \underline{s}_{2} \widehat{m}_{2} - \overline{s}_{1} \widehat{l}_{1} - \overline{s}_{2} \widehat{l}_{2}, y_{1} - \widehat{m}_{1} + \widehat{l}_{1}, y_{2} - \widehat{m}_{1} + \widehat{l}_{2}) = \nonumber \\
&& \varphi (x, y_{1}, y_{2}, \underline{s}_{1}, \underline{s}_{2}, \overline{s}_{1}, \overline{s}_{2}, \widehat{l}_{1}, \widehat{l}_{2}, \widehat{m}_{1}, \widehat{m}_{2}),
    \end{eqnarray*}
which completes the proof.
\end{proof}

\begin{lemma}\label{concave}
For a fixed $(x, y_{1}, y_{2}, \underline{s}_{1}, \underline{s}_{2}, \overline{s}_{1}, \overline{s}_{2}) \in \mathbb{R}_{+} \times \mathbb{R}_{+}^{2} \times \mathbb{D}$ the mapping
\begin{equation}
\mathbb{A}(x, y_{1}, y_{2}, \underline{s}_{1}, \underline{s}_{2}, \overline{s}_{1}, \overline{s}_{2})\ni (l_1,l_2,m_1,m_2)\mapsto \varphi (x, y_{1}, y_{2}, \underline{s}_{1}, \underline{s}_{2}, \overline{s}_{1}, \overline{s}_{2},l_1,l_2,m_1,m_2)
\end{equation}
is strictly concave within the class of controls $(l_1,l_2,m_1,m_2)$ such that $l_1m_1+l_2m_2=0$. In other words for $\lambda\in (0,1)$ and $(l_1^1,l_2^1,m_1^1,m_2^1), (l_1^2,l_2^2,m_1^2,m_2^2)\in \mathbb{A}(x, y_{1}, y_{2}, \underline{s}_{1}, \underline{s}_{2}, \overline{s}_{1}, \overline{s}_{2})$ such that $l_1^1m_1^1+l_2^1m_2^1=0$, $l_1^2m_1^2+l_2^2m_2^2=0$ and
$(l_1^1,l_2^1,m_1^1,m_2^1)\neq(l_1^2,l_2^2,m_1^2,m_2^2)$ we have
\begin{eqnarray}\label{conc}
&&\varphi(x, y_{1}, y_{2}, \underline{s}_{1}, \underline{s}_{2}, \overline{s}_{1}, \overline{s}_{2},\lambda(l_1^1,l_2^1,m_1^1,m_2^1)+(1-\lambda)(l_1^2,l_2^2,m_1^2,m_2^2))> \nonumber   \\
&&\lambda \varphi(x, y_{1}, y_{2}, \underline{s}_{1}, \underline{s}_{2}, \overline{s}_{1}, \overline{s}_{2},l_1^1,l_2^1,m_1^1,m_2^1) + \nonumber \\
&& (1-\lambda)\varphi(x, y_{1}, y_{2}, \underline{s}_{1}, \underline{s}_{2}, \overline{s}_{1}, \overline{s}_{2},l_1^2,l_2^2,m_1^2,m_2^2).
\end{eqnarray}
\end{lemma}
\begin{proof} By \eqref{vphi} inequality \eqref{conc} can be rewritten as
\begin{eqnarray}\label{con}
&& g \Big( x + \underline{s}_{1} (\lambda m_{1}^1+(1-\lambda) m_{1}^2) + \underline{s}_{2} (\lambda m_{2}^1 + (1-\lambda) m_{2}^2) - \overline{s}_{1} (\lambda l_{1}^1 + (1-\lambda) l_{1}^2) - \nonumber \\
&&\overline{s}_{2} (\lambda l_{2}^1 + (1-\lambda) l_{2}^2),
y_{1} - \lambda m_{1}^1 - (1-\lambda) m_{1}^2+ \lambda l_{1}^1 + \nonumber \\
&& (1-\lambda) l_{1}^2, y_{2} - \lambda m_{2}^1 - (1-\lambda) m_{2}^2+ \lambda l_{2}^1 + (1-\lambda) l_{2}^2 \Big) > \nonumber \\
&& \lambda g \Big( x + \underline{s}_{1} m_{1}^1 + \underline{s}_{2} m_{2}^1 - \overline{s}_{1} l_{1}^1 - \overline{s}_{2} l_{2}^1, y_{1} - m_{1}^1 + l_{1}^1, y_{2}  - m_{2}^1 + l_{2}^1 \Big)+ \nonumber \\
&&(1-\lambda) g \Big( x + \underline{s}_{1} m_{1}^2 + \underline{s}_{2} m_{2}^2 - \overline{s}_{1} l_{1}^2 - \overline{s}_{2} l_{2}^2, y_{1} - m_{1}^2 + l_{1}^2, y_{2}  - m_{2}^2 + l_{2}^2 \Big).
\end{eqnarray}
Since in each pair $(l_1^1,m_1^1)$, $(l_1^2,m_1^2)$, $(l_2^1,m_2^1)$, $(l_2^2,m_2^2)$ at least one term is equal to $0$ we easily have that \eqref{con} follows directly from strict concavity of the function $g$.
\end{proof}
\end{lemma}
\begin{proposition}\label{Prop2.3}
     Let $(x, y_{1}, y_{2}, \underline{s}_{1}, \underline{s}_{2}, \overline{s}_{1}, \overline{s}_{2}) \in \mathbb{R}_{+} \times \mathbb{R}_{+}^{2} \times \mathbb{D}$. Then there exists a unique
    \begin{equation*}
        (l_{1}, l_{2}, m_{1}, m_{2}) \in \mathbb{A} (x, y_{1}, y_{2}, \underline{s}_{1}, \underline{s}_{2}, \overline{s}_{1}, \overline{s}_{2})
    \end{equation*}
    such that
    \begin{equation*}
    \Phi (x, y_{1}, y_{2},  (l_{1}, l_{2}, m_{1}, m_{2})  = \varphi (x, y_{1}, y_{2}, \underline{s}_{1}, \underline{s}_{2}, \overline{s}_{1}, \overline{s}_{2}, l_{1}, l_{2}, m_{1}, m_{2}),
    \end{equation*}
    and  $l_{1} m_{1} + l_{2} m_{2} = 0$. If for $ (l_{1}, l_{2}, m_{1}, m_{2}) \in \mathbb{A} (x, y_{1}, y_{2}, \underline{s}_{1}, \underline{s}_{2}, \overline{s}_{1}, \overline{s}_{2})$ with either $\underline{s}_{1}< \overline{s}_{1}$ or $\underline{s}_{2}< \overline{s}_{2}$ supremum in $\Phi$ is attained then $l_{1} m_{1} + l_{2} m_{2} = 0$ .
\end{proposition}
\begin{proof} Existence follows directly from Lemma \ref{Lem2.01} and Lemma \ref {Lem2.1}.
     Uniqueness is a direct consequence of strict concavity of $g$ (see Lemma \ref{concave}). Note that when  $\underline{s}_{1}= \overline{s}_{1}$ or  $\underline{s}_{2}= \overline{s}_{2}$ then simultaneous buying and selling the same asset does not change the value function of $\phi$.
\end{proof}

Directly from Lemma \ref{Lem1.3} we obtain the following useful result
\begin{lemma}\label{usLem}
Assume $(x, y_{1}, y_{2}) \in \mathbb{R}_{+} \times \mathbb{R}_{+}^{2}$ and
    $(\underline{s}_{1}^{1}, \underline{s}_{2}^{1}, \overline{s}_{1}^{1}, \overline{s}_{2}^{1}) \in \mathbb{D}$, $ (\underline{s}_{1}^{2}, \underline{s}_{2}^{2}, \overline{s}_{1}^{2}, \overline{s}_{2}^{2}) \in \mathbb{D}$
are such that
    $\underline{s}_{1}^{1} \leq \underline{s}_{1}^{2} \leq \overline{s}_{1}^{2} \leq \overline{s}_{1}^{1}$ $\underline{s}_{2}^{1} \leq \underline{s}_{2}^{2} \leq \overline{s}_{2}^{2} \leq \overline{s}_{2}^{1}$.
Then
    $\Phi (x, y_{1}, y_{2}, \underline{s}_{1}^{1}, \underline{s}_{2}^{1}, \overline{s}_{1}^{1}, \overline{s}_{2}^{1}) \leq \Phi (x, y_{1}, y_{2}, \underline{s}_{1}^{2}, \underline{s}_{2}^{2}, \overline{s}_{1}^{2}, \overline{s}_{2}^{2})$ .
\end{lemma}

For every $(x, y_{1}, y_{2}, \underline{s}_{1}, \underline{s}_{2}, \overline{s}_{1}, \overline{s}_{2}) \in \mathbb{R}_{+} \times \mathbb{R}_{+}^{2} \times \mathbb{D}$ define the set
    \begin{eqnarray}
        && \overleftrightarrow{\mathbb{A}}(x, y_{1}, y_{2}, \underline{s}_{1}, \underline{s}_{2}, \overline{s}_{1}, \overline{s}_{2}) :=  \nonumber \\
        && \Big\{ (l_{1}, l_{2}, m_{1}, m_{2}) \in \mathbb{A}(x, y_{1}, y_{2}, \underline{s}_{1}, \underline{s}_{2}, \overline{s}_{1}, \overline{s}_{2}): \quad l_{1} m_{1} + l_{2} m_{2} = 0 \Big\} .
    \end{eqnarray}
We have

\begin{proposition}\label{continuity of widtilde}
There exists a function  $u : \mathbb{R}_{+} \times \mathbb{R}_{+}^{2} \times \mathbb{D} \rightarrow \mathbb{R}_{+} \times \mathbb{R}_{+}^{2}$ be such that for every $(x, y_{1}, y_{2}, \underline{s}_{1}, \underline{s}_{2}, \overline{s}_{1}, \overline{s}_{2}) \in \mathbb{R}_{+} \times \mathbb{R}_{+}^{2} \times \mathbb{D}$ we have that
    \begin{equation*}
        u(x, y_{1}, y_{2}, \underline{s}_{1}, \underline{s}_{2}, \overline{s}_{1}, \overline{s}_{2}) \in \overleftrightarrow{\mathbb{A}}(x, y_{1}, y_{2}, \underline{s}_{1}, \underline{s}_{2}, \overline{s}_{1}, \overline{s}_{2})
    \end{equation*}
    and
    \begin{equation*}
        \Phi (x, y_{1}, y_{2}, \underline{s}_{1}, \underline{s}_{2}, \overline{s}_{1}, \overline{s}_{2}) = \varphi \big(x, y_{1}, y_{2}, \underline{s}_{1}, \underline{s}_{2}, \overline{s}_{1}, \overline{s}_{2}, u(x, y_{1}, y_{2}, \underline{s}_{1}, \underline{s}_{2}, \overline{s}_{1}, \overline{s}_{2}) \big) .
    \end{equation*}
Furthermore $u$ is continuous and uniquely defined.
\end{proposition}
\begin{proof} Notice first that following the proof of Theorem 2.1 of \cite{RS1} we can show that the mapping $(x, y_{1}, y_{2}, \underline{s}_{1}, \underline{s}_{2}, \overline{s}_{1}, \overline{s}_{2}) \longmapsto \overleftrightarrow{\mathbb{A}}(x, y_{1}, y_{2}, \underline{s}_{1}, \underline{s}_{2}, \overline{s}_{1}, \overline{s}_{2})$ is continuous in Hausdorff metric. By Proposition \ref{Prop2.3} function $u$ is uniquely defined. Using Theorem \ref{select} we have continuity of $u$.
\end{proof}
Using Theorem \ref{select} one can show the following result

\begin{corollary}\label{Corollary on the continuity of the function widehat}
     For a fixed $(x, y_{1}, y_{2}, \underline{s}_{1}, \underline{s}_{2}, \overline{s}_{1}, \overline{s}_{2}) \in \mathbb{R}_{+} \times \mathbb{R}_{+}^{2} \times \mathbb{D}
    $ consider  $\widehat{u} : [\underline{s}_{1}, \overline{s}_{1}] \longrightarrow \mathbb{R}_{+}^{2} \times \mathbb{R}_{+}^{2}$ such that for every $s \in [\underline{s}_{1}, \overline{s}_{1}]$ we have
         $  \widehat{u}(s) \in \overleftrightarrow{\mathbb{A}}(x, y_{1}, y_{2}, s, \underline{s}_{2}, s, \overline{s}_{2})$  and
    \begin{equation*}
        \Phi (x, y_{1}, y_{2}, s, \underline{s}_{2}, s, \overline{s}_{2}) = \varphi \big(x, y_{1}, y_{2}, s, \underline{s}_{2}, s, \overline{s}_{2}, \widehat{u}(s) \big) .
    \end{equation*}
    Then $\widehat{u}$ is a  continuous function.
\end{corollary}

\section{Decomposition of the set of admissible investment strategies}

For every $(i_{1}, i_{2}) \in \{ -1, 0, 1 \}^{2}$ and for every $(x, y_{1}, y_{2}, \underline{s}_{1}, \underline{s}_{2}, \overline{s}_{1}, \overline{s}_{2}) \in \mathbb{R}_{+} \times \mathbb{R}_{+}^{2} \times \mathbb{D}$  define
\begin{eqnarray}
&& \mathbb{A}^{(i_{1}, i_{2})} (x, y_{1}, y_{2}, \underline{s}_{1}, \underline{s}_{2}, \overline{s}_{1}, \overline{s}_{2})
        := \Big\{ (l_{1}, l_{2}, m_{1}, m_{2}) \in \mathbb{A}(x, y_{1}, y_{2}, \underline{s}_{1}, \underline{s}_{2}, \overline{s}_{1}, \overline{s}_{2}) : \nonumber \\
&&\quad i_{1} = \sgn (l_{1} - m_{1}), i_{2} = \sgn (l_{2} - m_{2}), l_{1} m_{1} + l_{2} m_{2} = 0 \Big\}.
\end{eqnarray}
In particular this means that $\mathbb{A}^{(0, 1)}$ consists of the strategies no transaction with the first asset and buying the second one, while $\mathbb{A}^{(1, -1)}$
means that we buy the first asset and sell the second one.
\begin{lemma}\label{Lem3.1}
    Let $(x, y_{1}, y_{2}, \underline{s}_{1}, \underline{s}_{2}, \overline{s}_{1}, \overline{s}_{2}) \in \mathbb{R}_{+} \times \mathbb{R}_{+}^{2} \times \mathbb{D}$ and $(i_{1}, i_{2}) \in \{ - 1, 0, 1 \}^{2}$. Then the set $\mathbb{A}^{(i_{1}, i_{2})} (x, y_{1}, y_{2}, \underline{s}_{1}, \underline{s}_{2}, \overline{s}_{1}, \overline{s}_{2})$ is convex.
\end{lemma}

\begin{proof}
    Let $(l_{1}^{(1)}, l_{2}^{(1)}, m_{1}^{(1)}, m_{2}^{(1)}), (l_{1}^{(2)}, l_{2}^{(2)}, m_{1}^{(2)}, m_{2}^{(2)}) \in \mathbb{A}^{(i_{1}, i_{2})} (x, y_{1}, y_{2}, \underline{s}_{1}, \underline{s}_{2}, \overline{s}_{1}, \overline{s}_{2})$ and $\lambda \in [0, 1]$ be arbitrary.

    We will show that
    \begin{equation}\label{convexity of A^{(i_{1}, i_{2})}}
        \lambda \big( l_{1}^{(1)}, l_{2}^{(1)}, m_{1}^{(1)}, m_{2}^{(1)} \big) + (1 - \lambda ) \big( l_{1}^{(2)}, l_{2}^{(2)}, m_{1}^{(2)}, m_{2}^{(2)} \big) \in \mathbb{A}^{(i_{1}, i_{2})} (x, y_{1}, y_{2}, \underline{s}_{1}, \underline{s}_{2}, \overline{s}_{1}. \overline{s}_{2}) .
    \end{equation}

    It is rather clear that
    $i_{1} = \sgn \Big( \big[ \lambda l_{1}^{(1)} + (1 - \lambda ) l_{1}^{(2)} \big] - \big[ \lambda m_{1}^{(1)} + (1 - \lambda ) m_{1}^{(2)} \big] \Big)$  and  $i_{2} = \sgn \Big( \big[ \lambda l_{2}^{(1)} + (1 - \lambda ) l_{2}^{(2)} \big] - \big[ \lambda m_{2}^{(1)} + (1 - \lambda ) m_{2}^{(2)} \big] \Big)$ .
    Taking into account $\sgn(l_{1}^{(1)})=\sgn(l_{1}^{(2)})$, $\sgn(m_{1}^{(1)})=\sgn(m_{1}^{(2)})$,  $\sgn(l_{2}^{(1)})=\sgn(l_{2}^{(2)})$ and $\sgn(m_{2}^{(1)})=\sgn(m_{2}^{(2)})$ we obtain that
    \begin{equation*}
        \big[ \lambda l_{1}^{(1)} + (1 - \lambda ) l_{1}^{(2)} \big] \cdot \big[ \lambda m_{1}^{(1)} + (1 - \lambda ) m_{1}^{(2)} \big] + \big[ \lambda l_{2}^{(1)} + (1 - \lambda ) l_{2}^{(2)} \big] \cdot \big[ \lambda m_{2}^{(1)} + (1 - \lambda ) m_{2}^{(2)} \big] = 0 .
    \end{equation*}
which completes the proof of $(\ref{convexity of A^{(i_{1}, i_{2})}})$.
\end{proof}

\begin{remark}
The set $\mathbb{A}^{(i_{1}, i_{2})} (x, y_{1}, y_{2}, \underline{s}_{1}, \underline{s}_{2}, \overline{s}_{1}, \overline{s}_{2})$ is not necessarily closed. One can notice that the set
$\mathbb{A}^{(0, 0)} (x, y_{1}, y_{2}, \underline{s}_{1}, \underline{s}_{2}, \overline{s}_{1}, \overline{s}_{2})$ is closed and the set

\begin{eqnarray} \mathbb{A}^{(i_{1}, i_{2})} (x, y_{1}, y_{2}, \underline{s}_{1}, &&\underline{s}_{2}, \overline{s}_{1}, \overline{s}_{2})\cup \mathbb{A}^{(i_{1}, 0)} (x, y_{1}, y_{2}, \underline{s}_{1}, \underline{s}_{2}, \overline{s}_{1}, \overline{s}_{2})\cup \nonumber \\
&& \mathbb{A}^{(0, i_{2})} (x, y_{1}, y_{2}, \underline{s}_{1}, \underline{s}_{2}, \overline{s}_{1}, \overline{s}_{2})\cup \mathbb{A}^{(0, 0)} (x, y_{1}, y_{2}, \underline{s}_{1}, \underline{s}_{2}, \overline{s}_{1}, \overline{s}_{2})
\end{eqnarray}
 is closed as well.
\end{remark}

\begin{lemma}\label{Lem3.3}
    Let $(x, y_{1}, y_{2}, \underline{s}_{1}, \underline{s}_{2}, \overline{s}_{1}, \overline{s}_{2}) \in \mathbb{R}_{+} \times \mathbb{R}_{+}^{2} \times \mathbb{D}$ and $(i_{1}, i_{2}) \in \{ -1, 0, 1 \}^{2}$. Then the set
    \begin{equation*}
        \bigcup_{(j_{1}, j_{2}) \in \{ -1, 0, 1 \}^{2}, |j_{1}| + |j_{2}| \leq |i_{1}| + |i_{2}| } \mathbb{A}^{(j_{1}, j_{2})} (x, y_{1}, y_{2}, \underline{s}_{1}, \underline{s}_{2}, \overline{s}_{1}, \overline{s}_{2})
    \end{equation*}
    is closed and therefore compact.
\end{lemma}

\begin{proof}
    Let $(l_{1}^{n}, l_{2}^{n}, m_{1}^{n}, m_{2}^{n})$ be a sequence from the set defined in the above Lemma which is convergent to some
    \begin{equation*}
        (l_{1}, l_{2}, m_{1}, m_{2}) \in \mathbb{R}_{+}^{2} \times \mathbb{R}_{+}^{2}.
    \end{equation*}
    We show that $(l_{1}, l_{2}, m_{1}, m_{2})$ is also an element of this set.

     One can notice that there exists a subsequence $(n_{k})_{k = 1}^{\infty}$  such that each of the sequence $(l_{1}^{n_{k}} - m_{1}^{n_{k}})_{k = 1}^{\infty}$ and $(l_{2}^{n_{k}} - m_{2}^{n_{k}})_{k = 1}^{\infty}$ has the same sign.

    Define $r_{1} := \sgn (l_{1} - m_{1})$ and $r_{2} := \sgn (l_{2} - m_{2})$.

    Clearly  $|r_{1}| + |r_{2}| \leq |i_{1}| + |i_{2}|$.

     For each $n \in \mathbb{N}$ we have that $l_{1}^{n} m_{1}^{n} + l_{2}^{n} m_{2}^{n} = 0$ and consequently  $l_{1} m_{1} + l_{2} m_{2} = 0$. Therefore the considered set is closed.
\end{proof}

\begin{corollary}\label{Cor3.4}
     For each $(x, y_{1}, y_{2}, \underline{s}_{1}, \underline{s}_{2}, \overline{s}_{1}, \overline{s}_{2}) \in \mathbb{R}_{+} \times \mathbb{R}_{+}^{2} \times \mathbb{D}$ there exist a unique $(i_{1}, i_{2}) \in \{ - 1, 0, 1 \}^{2}$ and a unique $(l_{1}, l_{2}, m_{1}, m_{2}) \in \mathbb{A}^{(i_1,i_2)} (x, y_{1}, y_{2}, \underline{s}_{1}, \underline{s}_{2}, \overline{s}_{1}, \overline{s}_{2})$ such that
    \begin{equation*}
        \Phi (x, y_{1}, y_{2}, \underline{s}_{1}, \underline{s}_{2}, \overline{s}_{1}, \overline{s}_{2}) = \varphi (x, y_{1}, y_{2}, \underline{s}_{1}, \underline{s}_{2}, \overline{s}_{1}, \overline{s}_{2}, l_{1}, l_{2}, m_{1}, m_{2}) .
    \end{equation*}
\end{corollary}

\begin{proof}
    By Proposition $\ref{Prop2.3}$ there is a unique $(l_{1}, l_{2}, m_{1}, m_{2}) \in \mathbb{A} (x, y_{1}, y_{2}, \underline{s}_{1}, \underline{s}_{2}, \overline{s}_{1}, \overline{s}_{2})$ such that $l_1m_1+l_2m_2=0$, which is optimal. Therefore $i_1=\sgn(l_1-m_1)$ and $i_2=\sgn(l_2-m_2)$ is uniquely defined and  $(l_{1}, l_{2}, m_{1}, m_{2}) \in \mathbb{A}^{(i_1,i_2)} (x, y_{1}, y_{2}, \underline{s}_{1}, \underline{s}_{2}, \overline{s}_{1}, \overline{s}_{2})$.
\end{proof}

\section{Trading sets}

For every $(\underline{s}_{1}, \underline{s}_{2}, \overline{s}_{1}, \overline{s}_{2}) \in \mathbb{D}$ and for every $(i_{1}, i_{2}) \in \{ - 1, 0, 1 \}^{2}$ define the following trading set
\begin{eqnarray}
    &&\mathbf{t}^{(i_{1}, i_{2})}  (\underline{s}_{1}, \underline{s}_{2}, \overline{s}_{1}, \overline{s}_{2}) := \Big\{ (x, y_{1}, y_{2}) \in \mathbb{R}_{+} \times \mathbb{R}_{+}^{2} : \nonumber \\
    &&\exists_{(l_{1}, l_{2}, m_{1}, m_{2}) \in \mathbb{A}^{(i_{1}, i_{2})} (x, y_{1}, y_{2}, \underline{s}_{1}, \underline{s}_{2}, \overline{s}_{1}, \overline{s}_{2})} \Phi (x, y_{1}, y_{2}, \underline{s}_{1}, \underline{s}_{2}, \overline{s}_{1}, \overline{s}_{2}) =  \nonumber \\
    && \varphi (x, y_{1}, y_{2}, \underline{s}_{1}, \underline{s}_{2}, \overline{s}_{1}, \overline{s}_{2}, l_{1}, l_{2}, m_{1}, m_{2}) \Big\}.
\end{eqnarray}

\begin{lemma}\label{disjoint}
    Let $(\underline{s}_{1}, \underline{s}_{2}, \overline{s}_{1}, \overline{s}_{2}) \in \mathbb{D}$ and $(i_{1}, i_{2}) \in \{ - 1, 0, 1 \}^{2}$. Then $\mathbf{t}^{(i_{1}, i_{2})} (\underline{s}_{1}, \underline{s}_{2}, \overline{s}_{1}, \overline{s}_{2})$ is a Borel subset of $ \mathbb{R}_{+} \times \mathbb{R}_{+}^{2}$, which we denote that it is an element of $\mathcal{B} \big( \mathbb{R}_{+} \times \mathbb{R}_{+}^{2} \big)$. Furthermore $\Big\{ \mathbf{t}^{(i_{1}, i_{2})} (\underline{s}_{1}, \underline{s}_{2}, \overline{s}_{1}, \overline{s}_{2}) \Big\}_{(i_{1}, i_{2}) \in \{ - 1, 0, 1 \}}$ is a family of disjoint sets which cover $\mathbb{R}_{+} \times \mathbb{R}_{+}^{2}$.
\end{lemma}
\begin{proof}
For measurability of $\mathbf{t}^{(i_{1}, i_{2})} (\underline{s}_{1}, \underline{s}_{2}, \overline{s}_{1}, \overline{s}_{2})$ it suffices to show that
    \begin{equation*}
        \Big( \mathbb{R}_{+} \times \mathbb{R}_{+}^{2} \Big) \setminus \mathbf{t}^{(i_{1}, i_{2})} (\underline{s}_{1}, \underline{s}_{2}, \overline{s}_{1}, \overline{s}_{2}) \in \mathcal{B} \big( \mathbb{R}_{+} \times \mathbb{R}_{+}^{2} \big) .
    \end{equation*}
     We have
    \begin{eqnarray*}
&&  \Big( \mathbb{R}_{+}  \times \mathbb{R}_{+}^{2} \Big) \setminus \mathbf{t}^{(i_{1}, i_{2})} (\underline{s}_{1}, \underline{s}_{2}, \overline{s}_{1}, \overline{s}_{2}) =  \nonumber \\
 &&\Big\{ (x, y_{1}, y_{2}) \in \mathbb{R}_{+} \times \mathbb{R}_{+}^{2} : \forall_{(l_{1}, l_{2}, m_{1}, m_{2}) \in \mathbb{A}^{(i_{1}, i_{2})} (x, y_{1}, y_{2}, \underline{s}_{1}, \underline{s}_{2}, \overline{s}_{1}, \overline{s}_{2})}
         \nonumber \\
          &&\Phi (x, y_{1}, y_{2}, \underline{s}_{1}, \underline{s}_{2}, \overline{s}_{1}, \overline{s}_{2}) > \varphi (x, y_{1}, y_{2}, \underline{s}_{1}, \underline{s}_{2}, \overline{s}_{1}, \overline{s}_{2}, l_{1}, l_{2}, m_{1}, m_{2}) \Big\} = \\        
       && =  \Big\{ (x, y_{1}, y_{2}) \in \mathbb{R}_{+} \times \mathbb{R}_{+}^{2} :
         \forall_{(l_{1}, l_{2}, m_{1}, m_{2}) \in \mathbb{A}^{(i_{1}, i_{2})} (x, y_{1}, y_{2}, \underline{s}_{1}, \underline{s}_{2}, \overline{s}_{1}, \overline{s}_{2}) \cap \mathbb{Q}^{4}}  \nonumber \\
         && \Phi (x, y_{1}, y_{2}, \underline{s}_{1}, \underline{s}_{2}, \overline{s}_{1}, \overline{s}_{2}) > \varphi (x, y_{1}, y_{2}, \underline{s}_{1}, \underline{s}_{2}, \overline{s}_{1}, \overline{s}_{2}, l_{1}, l_{2}, m_{1}, m_{2}) \Big\} = \\
       && =  \bigcap_{(l_{1}, l_{2}, m_{1}, m_{2}) \in \mathbb{A}^{(i_{1}, i_{2})} (x, y_{1}, y_{2}, \underline{s}_{1}, \underline{s}_{2}, \overline{s}_{1}, \overline{s}_{2}) \cap \mathbb{Q}^{4}} \Big\{ (x, y_{1}, y_{2}) \in \mathbb{R}_{+} \times \mathbb{R}_{+}^{2} : \\
 &&  \Phi (x, y_{1}, y_{2}, \underline{s}_{1}, \underline{s}_{2}, \overline{s}_{1}, \overline{s}_{2}) > \varphi (x, y_{1}, y_{2}, \underline{s}_{1}, \underline{s}_{2}, \overline{s}_{1}, \overline{s}_{2}, l_{1}, l_{2}, m_{1}, m_{2}) \Big\} = \\
       && =  \bigcap_{(l_{1}, l_{2}, m_{1}, m_{2}) \in \mathbb{A}^{(i_{1}, i_{2})} (x, y_{1}, y_{2}, \underline{s}_{1}, \underline{s}_{2}, \overline{s}_{1}, \overline{s}_{2}) \cap \mathbb{Q}^{4}} \quad \bigcup_{q \in \mathbb{Q}, q < \Phi (x, y_{1}, y_{2}, \underline{s}_{1}, \underline{s}_{2}, \overline{s}_{1}, \overline{s}_{2})} \Big\{ (x, y_{1}, y_{2}) \in \mathbb{R}_{+} \times \mathbb{R}_{+}^{2} : \\
        &&  q > \varphi (x, y_{1}, y_{2}, \underline{s}_{1}, \underline{s}_{2}, \overline{s}_{1}, \overline{s}_{2}, l_{1}, l_{2}, m_{1}, m_{2}) \Big\} .
    \end{eqnarray*}
where $\mathbb{Q}$ stands for rationals. Consequently the above set is in  the $\sigma$-field $\mathcal{B} \big( \mathbb{R}_{+} \times \mathbb{R}_{+}^{2} \big)$ and the same concerns the set $\mathbf{t}^{(i_{1}, i_{2})} (\underline{s}_{1}, \underline{s}_{2}, \overline{s}_{1}, \overline{s}_{2})$, which completes the proof of measurability. The fact that the family $\Big\{ \mathbf{t}^{(i_{1}, i_{2})} (\underline{s}_{1}, \underline{s}_{2}, \overline{s}_{1}, \overline{s}_{2}) \Big\}_{(i_{1}, i_{2}) \in \{ - 1, 0, 1 \}}$ consists of disjoint sets follows directly from Corollary $\ref{Cor3.4}$.
\end{proof}
\begin{proposition}\label{top1-3}
For $(\underline{s}_{1}, \underline{s}_{2}, \overline{s}_{1}, \overline{s}_{2}) \in \mathbb{D}$ and $i,j \in \{ - 1, 0, 1 \}$ we have the following properties of trading zones:
    \begin{equation}\label{prop1}
         \mathbf{t}^{(-1, i)}(\underline{s}_{1}, \underline{s}_{2}, \underline{s}_{1}, \overline{s}_{2}) \subseteq \mathbf{t}^{(-1, i)}(\underline{s}_{1}, \underline{s}_{2}, \overline{s}_{1}, \overline{s}_{2}),
    \end{equation}
\begin{equation}\label{prop2}
  \mathbf{t}^{(0, i)}(\underline{s}_{1}, \underline{s}_{2}, \underline{s}_{1}, \overline{s}_{2}) \subseteq \mathbf{t}^{(0, i)}(\underline{s}_{1}, \underline{s}_{2}, \overline{s}_{1}, \overline{s}_{2}),
\end{equation}
\begin{equation}\label{prop3}
 \mathbf{t}^{(0, i)}(\underline{s}_{1}, \underline{s}_{2}, \underline{s}_{1}, \overline{s}_{2}) \cap \mathbf{t}^{(-1, j)}(\underline{s}_{1}, \underline{s}_{2}, \overline{s}_{1}, \overline{s}_{2})=\emptyset.
  \end{equation}
 \end{proposition}
\begin{proof}
Assume $(x,y_1,y_2)\in \mathbf{t}^{(-1, i)}(\underline{s}_{1}, \underline{s}_{2}, \underline{s}_{1}, \overline{s}_{2})$.

Then there is $(0,l_2,m_1,m_2)\in$ $\mathbb{A}(x, y_{1}, y_{2}, \underline{s}_{1}, \underline{s}_{2}, \underline{s}_{1}, \overline{s}_{2})$ such that
\begin{equation*}
 \Phi (x, y_{1}, y_{2}, \underline{s}_{1}, \underline{s}_{2}, \underline{s}_{1}, \overline{s}_{2}) = \varphi \big( x, y_{1}, y_{2}, \underline{s}_{1}, \underline{s}_{2}, \underline{s}_{1}, \overline{s}_{2}, 0, l_2, m_1, m_2 \big)
\end{equation*}
We also have that $(0,l_2,m_1,m_2)\in \mathbb{A}(x, y_{1}, y_{2}, \underline{s}_{1}, \underline{s}_{2}, \overline{s}_{1}, \overline{s}_{2})$ and consequently
\begin{eqnarray*}
&&\Phi (x, y_{1}, y_{2}, \underline{s}_{1}, \underline{s}_{2}, \overline{s}_{1}, \overline{s}_{2})\leq
\Phi (x, y_{1}, y_{2}, \underline{s}_{1}, \underline{s}_{2}, \underline{s}_{1}, \overline{s}_{2})= \\
&&\varphi \big( x, y_{1}, y_{2}, \underline{s}_{1}, \underline{s}_{2}, \underline{s}_{1}, \overline{s}_{2}, 0, l_2, m_1, m_2 \big)=\varphi \big( x, y_{1}, y_{2}, \underline{s}_{1}, \underline{s}_{2}, \overline{s}_{1}, \overline{s}_{2}, 0, l_2, m_1, m_2 \big)\leq \\
&& \Phi (x, y_{1}, y_{2}, \underline{s}_{1}, \underline{s}_{2}, \overline{s}_{1}, \overline{s}_{2})
\end{eqnarray*}
which means that $\Phi (x, y_{1}, y_{2}, \underline{s}_{1}, \underline{s}_{2}, \overline{s}_{1}, \overline{s}_{2})=\varphi \big( x, y_{1}, y_{2}, \underline{s}_{1}, \overline{s}_{2}, \underline{s}_{1}, \overline{s}_{2}, 0, l_2, m_1, m_2 \big)$ and
\begin{equation*}(x,y_1,y_2)\in \mathbf{t}^{(-1, i)}(\underline{s}_{1}, \underline{s}_{2}, \overline{s}_{1}, \overline{s}_{2}).
\end{equation*}
To show \eqref{prop2} assume that $(x,y_1,y_2)\in \mathbf{t}^{( 0, i)}(\underline{s}_{1}, \underline{s}_{2}, \underline{s}_{1}, \overline{s}_{2})$ and there exists
\begin{equation*}(0,l_2,0,m_2)\in \mathbb{A}(x, y_{1}, y_{2}, \underline{s}_{1}, \underline{s}_{2}, \underline{s}_{1}, \overline{s}_{2})
\end{equation*}
 such that
\begin{equation*}
 \Phi (x, y_{1}, y_{2}, \underline{s}_{1}, \underline{s}_{2}, \underline{s}_{1}, \overline{s}_{2}) = \varphi \big( x, y_{1}, y_{2}, \underline{s}_{1}, \underline{s}_{2}, \underline{s}_{1}, \overline{s}_{2}, 0, l_2, 0, m_2 \big)
\end{equation*}
Clearly $(0,l_2,0,m_2)\in \mathbb{A}(x, y_{1}, y_{2}, \underline{s}_{1}, \underline{s}_{2}, \overline{s}_{1}, \overline{s}_{2})$ and
\begin{eqnarray*}
&&\Phi (x, y_{1}, y_{2}, \underline{s}_{1}, \underline{s}_{2}, \overline{s}_{1}, \overline{s}_{2})\leq
\Phi (x, y_{1}, y_{2}, \underline{s}_{1}, \underline{s}_{2}, \underline{s}_{1}, \overline{s}_{2})= \\
&&\varphi \big( x, y_{1}, y_{2}, \underline{s}_{1}, \underline{s}_{2}, \underline{s}_{1}, \overline{s}_{2}, 0, l_2, 0, m_2 \big)=\varphi \big( x, y_{1}, y_{2}, \underline{s}_{1}, \underline{s}_{2}, \overline{s}_{1}, \overline{s}_{2}, 0, l_2, 0, m_2 \big)\leq \\
&&\Phi (x, y_{1}, y_{2}, \underline{s}_{1}, \underline{s}_{2}, \overline{s}_{1}, \overline{s}_{2}).
\end{eqnarray*}
Therefore $\Phi (x, y_{1}, y_{2}, \underline{s}_{1}, \underline{s}_{2}, \overline{s}_{1}, \overline{s}_{2})=\varphi \big( x, y_{1}, y_{2}, \underline{s}_{1}, \overline{s}_{2}, \underline{s}_{1}, \overline{s}_{2}, 0, l_2, 0, m_2 \big)$ and
\begin{equation*}(x,y_1,y_2)\in \mathbf{t}^{(0, i)}(\underline{s}_{1}, \underline{s}_{2}, \overline{s}_{1}, \overline{s}_{2}).
\end{equation*}
To show \eqref{prop3} notice that by \eqref{prop2}
$ \mathbf{t}^{(0, i)}(\underline{s}_{1}, \underline{s}_{2}, \underline{s}_{1}, \overline{s}_{2})\subset \mathbf{t}^{(0, i)}(\underline{s}_{1}, \underline{s}_{2}, \overline{s}_{1}, \overline{s}_{2})$
 and by Lemma \ref{disjoint} the sets $\mathbf{t}^{(-1, j)}(\underline{s}_{1}, \underline{s}_{2}, \overline{s}_{1}, \overline{s}_{2})$ and $\mathbf{t}^{(0, i)}(\underline{s}_{1}, \underline{s}_{2}, \overline{s}_{1}, \overline{s}_{2})$ are disjoint.
\end{proof}
\begin{remark}
The first two properties of Proposition \ref{top1-3} are intuitive. In the case of \eqref{prop1} when it is optimal to sell the first asset in the case when the bid and ask price are equal to $ \underline{s}_{1}$, then it is also optimal to sell it when the ask price increases to $\overline{s}_{1}$. In the case of \eqref{prop2} when it is optimal to do nothing with the first asset in the case when the bid and ask price is $ \underline{s}_{1}$, then it is also optimal to do nothing with the first asset when the ask price increases to $\overline{s}_{1}$.
\end{remark}

\begin{corollary}
For $(\underline{s}_{1}, \underline{s}_{2}, \overline{s}_{1}, \overline{s}_{2}) \in \mathbb{D}$ and $i\neq j \in \{ - 1, 0, 1 \}$
\begin{equation}\label{prop6}
\mathbf{t}^{(- 1, j)} (\underline{s}_{1}, \underline{s}_{2}, \underline{s}_{1}, \overline{s}_{2}) \cap \mathbf{t}^{(- 1, i)} (\underline{s}_{1}, \underline{s}_{2}, \overline{s}_{1}, \overline{s}_{2})=\emptyset.
    \end{equation}
\end{corollary}
\begin{proof}
Notice that by \eqref{prop1}
\begin{equation*}
\mathbf{t}^{(- 1, j)} (\underline{s}_{1}, \underline{s}_{2}, \underline{s}_{1}, \overline{s}_{2}) \subseteq \mathbf{t}^{(- 1, j)} (\underline{s}_{1}, \underline{s}_{2}, \overline{s}_{1}, \overline{s}_{2})
\end{equation*}
and by Lemma \ref{disjoint} $\Big\{ \mathbf{t}^{(i_{1}, i_{2})} (\underline{s}_{1}, \underline{s}_{2}, \overline{s}_{1}, \overline{s}_{2}) \Big\}_{(i_{1}, i_{2}) \in \{ - 1, 0, 1 \}}$ is a family of disjoint sets.
\end{proof}

\begin{proposition}\label{top4}
    Let $(\underline{s}_{1}, \underline{s}_{2}, \overline{s}_{1}, \overline{s}_{2}) \in \mathbb{D}$ and $i,j \in \{ - 1, 0, 1 \}$. Then
    \begin{equation}\label{prop4}
         \mathbf{t}^{( - 1, i)}(\underline{s}_{1}, \underline{s}_{2}, \overline{s}_{1}, \overline{s}_{2}) \cap \mathbf{t}^{(1, j)}(\underline{s}_{1}, \underline{s}_{2}, \underline{s}_{1}, \overline{s}_{2})= \emptyset.
    \end{equation}
\end{proposition}

\begin{proof}
Assume that $(x, y_{1}, y_{2}) \in \mathbf{t}^{( - 1, i)}(\underline{s}_{1}, \underline{s}_{2}, \overline{s}_{1}, \overline{s}_{2}) \cap \mathbf{t}^{(1, j)}(\underline{s}_{1}, \underline{s}_{2}, \underline{s}_{1}, \overline{s}_{2})$.
There exist
    $\big( 0, l_{2}^{(- 1)}, m_{1}^{(- 1)}, m_{2}^{(- 1)} \big) \in \mathbb{A}^{( - 1, i)}(x, y_{1}, y_{2}, \underline{s}_{1}, \underline{s}_{2}, \overline{s}_{1}, \overline{s}_{2} )$ and
     $\big( l_{1}^{(1)}, l_{2}^{(1)}, 0, m_{2}^{(1)} \big) \in \mathbb{A}^{( 1, j)}(x, y_{1}, y_{2},$ $\underline{s}_{1}, \underline{s}_{2}, \underline{s}_{1}, \overline{s}_{2})$
such that
\begin{equation*}
 \Phi (x, y_{1}, y_{2}, \underline{s}_{1}, \underline{s}_{2}, \overline{s}_{1}, \overline{s}_{2}) = \varphi \big( x, y_{1}, y_{2}, \underline{s}_{1}, \overline{s}_{2}, \overline{s}_{1}, \overline{s}_{2}, 0, l_{2}^{(- 1)}, m_{1}^{(- 1)}, m_{2}^{(- 1)} \big)
\end{equation*}
and
\begin{equation*}
\Phi (x, y_{1}, y_{2}, \underline{s}_{1}, \underline{s}_{2}, \underline{s}_{1}, \overline{s}_{2}) =
\varphi \big( x, y_{1}, y_{2}, \underline{s}_{1}, \underline{s}_{2}, \underline{s}_{1}, \overline{s}_{2}, l_{1}^{(1)}, l_{2}^{(1)}, 0, m_{2}^{(1)} \big).
\end{equation*}
Since  $\mathbb{A}(x, y_{1}, y_{2}, \underline{s}_{1}, \underline{s}_{2}, \overline{s}_{1}, \overline{s}_{2} )\subset \mathbb{A}(x, y_{1}, y_{2}, \underline{s}_{1}, \underline{s}_{2}, \underline{s}_{1}, \overline{s}_{2} )$ we have that
\begin{equation*}
\Phi (x, y_{1}, y_{2}, \underline{s}_{1}, \underline{s}_{2}, \underline{s}_{1}, \overline{s}_{2})\geq \Phi (x, y_{1}, y_{2}, \underline{s}_{1}, \underline{s}_{2}, \overline{s}_{1}, \overline{s}_{2}).
\end{equation*}
Moreover $m_{1}^{(- 1)} > 0 $ and $l_{1}^{(1)} > 0$.
Let
\begin{equation*}
    \lambda := \frac{l_{1}^{(1)}}{m_{1}^{(- 1)} + l_{1}^{(1)}} .
\end{equation*}
By convexity of the set $\mathbb{A}(x, y_{1}, y_{2}, \underline{s}_{1}, \underline{s}_{2}, \underline{s}_{1}, \overline{s}_{2} )$
\begin{eqnarray*}
  &&  (l_{1}, l_{2}, m_{1}, m_{2}) := \lambda \big( 0, l_{2}^{(- 1)}, m_{1}^{(- 1)}, m_{2}^{(- 1)} \big) + (1 - \lambda ) \big( l_{1}^{(1)}, l_{2}^{(1)}, 0, m_{2}^{(1)} \big) \\
  &&\in \mathbb{A}(x, y_{1}, y_{2}, \underline{s}_{1}, \underline{s}_{2}, \underline{s}_{1}, \overline{s}_{2} ).
\end{eqnarray*}
Since $(1-\lambda)l_{1}^{(1)}=\lambda m_{1}^{(- 1)}$ using Lemma \ref{Lem1.4} we have
$(0, l_{2}, 0, m_{2})\in \mathbb{A}(x, y_{1}, y_{2}, \underline{s}_{1}, \underline{s}_{2}, \underline{s}_{1}, \overline{s}_{2} )$. But then also $(0, l_{2}, 0, m_{2})\in \mathbb{A}(x, y_{1}, y_{2}, \underline{s}_{1}, \underline{s}_{2}, \overline{s}_{1}, \overline{s}_{2} )$.
Therefore by strict concavity of $g$, using Lemma \ref{concave} we have
\begin{eqnarray*}
&&\varphi \big( x, y_{1}, y_{2}, \underline{s}_{1}, \underline{s}_{2}, \underline{s}_{1}, \overline{s}_{2},0, l_{2},0, m_{2} \big)=\varphi \big( x, y_{1}, y_{2}, \underline{s}_{1}, \underline{s}_{2}, \underline{s}_{1}, \overline{s}_{2}, l_{1}, l_{2},m_{1}, m_{2} \big)> \\
&&\lambda \varphi \big( x, y_{1}, y_{2}, \underline{s}_{1}, \underline{s}_{2}, \underline{s}_{1}, \overline{s}_{2}, 0, l_{2}^{(- 1)}, m_{1}^{(- 1)}, m_{2}^{(- 1)} \big)+\\
&&(1-\lambda)\varphi \big( x, y_{1}, y_{2}, \underline{s}_{1}, \underline{s}_{2}, \underline{s}_{1}, \overline{s}_{2},l_{1}^{(1)}, l_{2}^{(1)}, 0, m_{2}^{(1)} \big)\geq
\Phi (x, y_{1}, y_{2}, \underline{s}_{1}, \underline{s}_{2}, \overline{s}_{1}, \overline{s}_{2}),
\end{eqnarray*}
which contradicts optimality of $\big( 0, l_{2}^{(- 1)}, m_{1}^{(- 1)}, m_{2}^{(- 1)} \big)$.
\end{proof}

\section{Equalities of selling and buying zones}
In this section we show several results concerning selling and buying zones. Main result of this section form Corollaries \ref{top8}, \ref{top9} and \ref{tobp8}, \ref{tobp9}. They are not used to construct shadow price, however they provide us with an interesting information about the properties of selling and buying zones. 
\begin{proposition}\label{top7}
For $(\underline{s}_{1}, \underline{s}_{2}, \overline{s}_{1}, \overline{s}_{2}) \in \mathbb{D}$ and $i \in \{ - 1, 0, 1 \}$
we have
 \begin{equation}\label{prop7}
        \mathbf{t}^{(- 1, i)} (\underline{s}_{1}, \underline{s}_{2}, \overline{s}_{1}, \overline{s}_{2}) \subseteq \mathbf{t}^{(- 1, i)} (\underline{s}_{1}, \underline{s}_{2}, \underline{s}_{1}, \overline{s}_{2}) .
    \end{equation}
\end{proposition}
\begin{proof}
Assume that $(x,y_1,y_2)\in \mathbf{t}^{(- 1, i)} (\underline{s}_{1}, \underline{s}_{2}, \overline{s}_{1}, \overline{s}_{2})$
and $(x,y_1,y_2)\notin \mathbf{t}^{(- 1, i)} (\underline{s}_{1}, \underline{s}_{2}, \underline{s}_{1}, \overline{s}_{2})$.
Then
\begin{eqnarray}\label{c1}
&& (x, y_{1}, y_{2})  \in \bigcup_{(k_{1}, k_{2}) \in \{ - 1, 0, 1 \}^{2} \setminus \{ (- 1, i) \}} \mathbf{t}^{(k_{1}, k_{2})} (\underline{s}_{1}, \underline{s}_{2}, \underline{s}_{1}, \overline{s}_{2}) =  \nonumber \\
&&\bigcup_{j \in \{ - 1, 0, 1 \} \setminus \{ i \}} \mathbf{t}^{(- 1, j)} (\underline{s}_{1}, \underline{s}_{2}, \underline{s}_{1}, \overline{s}_{2}) \cup \nonumber \\
&&  \bigcup_{j \in \{ - 1, 0, 1 \}} \mathbf{t}^{(0, j)} (\underline{s}_{1}, \underline{s}_{2}, \underline{s}_{1}, \overline{s}_{2}) \cup \bigcup_{j \in \{ - 1, 0, 1 \}} \mathbf{t}^{(1 , j)} (\underline{s}_{1}, \underline{s}_{2}, \underline{s}_{1}, \overline{s}_{2}) .
\end{eqnarray}
From \eqref{prop6}  we have that
    \begin{equation}\label{c2}
        (x, y_{1}, y_{2}) \notin \bigcup_{j \in \{ - 1, 0, 1 \} \setminus \{ i \}} \mathbf{t}^{(- 1, j)} (\underline{s}_{1}, \underline{s}_{2}, \underline{s}_{1}, \overline{s}_{2}) .
    \end{equation}
 From \eqref{prop3} we have that
    \begin{equation}\label{c3}
        (x, y_{1}, y_{2}) \notin \bigcup_{j \in \{ - 1, 0, 1 \}} \mathbf{t}^{(0, j)} (\underline{s}_{1}, \underline{s}_{2}, \underline{s}_{1}, \overline{s}_{2}) .
    \end{equation}
    and finally from \eqref{prop4}
     we have that
    \begin{equation}\label{c4}
        (x, y_{1}, y_{2}) \notin \bigcup_{j \in \{ - 1, 0, 1 \}} \mathbf{t}^{(1, j)} (\underline{s}_{1}, \underline{s}_{2}, \underline{s}_{1}, \overline{s}_{2}) .
    \end{equation}
Summarizing \eqref{c2}-\eqref{c4} we obtain a contradiction to \eqref{c1} and \eqref{prop7} is satisfied.
\end{proof}
Combining \eqref{prop1} and Proposition \ref{top7} we obtain
\begin{corollary}\label{top8}
For $(\underline{s}_{1}, \underline{s}_{2}, \overline{s}_{1}, \overline{s}_{2}) \in \mathbb{D}$ and $i \in \{ - 1, 0, 1 \}$ we
have
    \begin{equation}\label{prop8}
        \mathbf{t}^{(-1, i)} (\underline{s}_{1}, \underline{s}_{2}, \overline{s}_{1}, \overline{s}_{2}) = \mathbf{t}^{(-1, i)} (\underline{s}_{1}, \underline{s}_{2}, \underline{s}_{1}, \overline{s}_{2}) .
    \end{equation}
\end{corollary}
By analogy we can also obtain
\begin{corollary}\label{top9}
For $(\underline{s}_{1}, \underline{s}_{2}, \overline{s}_{1}, \overline{s}_{2}) \in \mathbb{D}$ and $i \in \{ - 1, 0, 1 \}$ we
have
    \begin{equation}\label{prop9}
        \mathbf{t}^{(i, -1)} (\underline{s}_{1}, \underline{s}_{2}, \overline{s}_{1}, \overline{s}_{2}) = \mathbf{t}^{(i, -1)} (\underline{s}_{1}, \underline{s}_{2}, \overline{s}_{1}, \underline{s}_{2}) .
    \end{equation}
\end{corollary}

Now we present several properties concerning mainly buying zones

\begin{proposition}\label{tobp1}
For $(\underline{s}_{1}, \underline{s}_{2}, \overline{s}_{1}, \overline{s}_{2}) \in \mathbb{D}$ and $i\neq j \in \{ - 1, 0, 1 \}$ we have

 \begin{equation}\label{bp1}
        \mathbf{t}^{(1, i)} (\overline{s}_{1}, \underline{s}_{2}, \overline{s}_{1}, \overline{s}_{2}) \subseteq \mathbf{t}^{(1, i)} (\underline{s}_{1}, \underline{s}_{2}, \overline{s}_{1}, \overline{s}_{2}),
    \end{equation}
    \begin{equation}\label{bp2}
        \mathbf{t}^{(0, i)} (\overline{s}_{1}, \underline{s}_{2}, \overline{s}_{1}, \overline{s}_{2}) \subseteq \mathbf{t}^{(0, i)} (\underline{s}_{1}, \underline{s}_{2}, \overline{s}_{1}, \overline{s}_{2}),
    \end{equation}
     \begin{equation}\label{bp3}
           \mathbf{t}^{(0, i)} (\overline{s}_{1}, \underline{s}_{2}, \overline{s}_{1}, \overline{s}_{2}) \cap \mathbf{t}^{(1, j)} (\underline{s}_{1}, \underline{s}_{2}, \overline{s}_{1}, \overline{s}_{2}) = \emptyset .
        \end{equation}
\end{proposition}
\begin{proof}
We use similar arguments as in the proof of Proposition \ref{top1-3}. Notice that it follows our intuition:  when it is optimal to  buy for $\overline{s}_1$ then it will be also optimal to buy for  $\underline{s}_1$. When we do nothing for bid and ask prices equal to $\overline{s}_1$ then we do not change the strategy when bid price decreases to $\underline{s}_1$. To show \eqref{bp3} notice that by \eqref{bp1}
 \begin{equation*}
 \mathbf{t}^{(0, i)} (\overline{s}_{1}, \underline{s}_{2}, \overline{s}_{1}, \overline{s}_{2})\subseteq \mathbf{t}^{(0, i)} (\underline{s}_{1}, \underline{s}_{2}, \overline{s}_{1}, \overline{s}_{2})
 \end{equation*}
 and use Lemma \ref{disjoint}.
\end{proof}

\begin{corollary}\label{tobp6}
For $(\underline{s}_{1}, \underline{s}_{2}, \overline{s}_{1}, \overline{s}_{2}) \in \mathbb{D}$ and $i\neq j \in \{ - 1, 0, 1 \}$
\begin{equation}\label{bp6}
\mathbf{t}^{(1, j)} (\underline{s}_{1}, \underline{s}_{2}, \overline{s}_{1}, \overline{s}_{2}) \cap \mathbf{t}^{(1, i)} (\overline{s}_{1}, \underline{s}_{2}, \overline{s}_{1}, \overline{s}_{2})=\emptyset.
    \end{equation}
\end{corollary}
\begin{proof}
Notice that by \eqref{bp1}
\begin{equation*}
\mathbf{t}^{(1, i)} (\overline{s}_{1}, \underline{s}_{2}, \overline{s}_{1}, \overline{s}_{2}) \subseteq \mathbf{t}^{(1, i)} (\underline{s}_{1}, \underline{s}_{2}, \overline{s}_{1}, \overline{s}_{2})
\end{equation*}
and by Lemma \ref{disjoint} $\Big\{ \mathbf{t}^{(i_{1}, i_{2})} (\underline{s}_{1}, \underline{s}_{2}, \overline{s}_{1}, \overline{s}_{2}) \Big\}_{(i_{1}, i_{2}) \in \{ - 1, 0, 1 \}}$ is a family of disjoint sets.
\end{proof}

 \begin{proposition}\label{tobp4}
For $(\underline{s}_{1}, \underline{s}_{2}, \overline{s}_{1}, \overline{s}_{2}) \in \mathbb{D}$ and $i, j \in \{ - 1, 0, 1 \}$ we have
        \begin{equation}\label{bp4}
            \mathbf{t}^{(1, i)} (\underline{s}_{1}, \underline{s}_{2}, \overline{s}_{1}, \overline{s}_{2}) \cap \mathbf{t}^{(- 1, j)} (\overline{s}_{1}, \underline{s}_{2}, \overline{s}_{1}, \overline{s}_{2}) = \emptyset .
        \end{equation}
    \end{proposition}
\begin{proof}
Assume there exists $(x, y_{1}, y_{2}) \in \mathbf{t}^{(1, i)} (\underline{s}_{1}, \underline{s}_{2}, \overline{s}_{1}, \overline{s}_{2}) \cap \mathbf{t}^{(- 1, j)} (\overline{s}_{1}, \underline{s}_{2}, \overline{s}_{1}, \overline{s}_{2})$.
Then there exists
        \begin{equation*}
            (\widehat{l}_{1}, \widehat{l}_{2}, 0, \widehat{m}_{2}) \in \mathbb{A}^{(1, i)}(x, y_{1}, y_{2}, \underline{s}_{1}, \underline{s}_{2}, \overline{s}_{1}, \overline{s}_{2})
        \end{equation*}
        such that
        \begin{equation*}
            \Phi (x, y_{1}, y_{2}, \underline{s}_{1}, \underline{s}_{2}, \overline{s}_{1}, \overline{s}_{2}) = \varphi (x, y_{1}, y_{2}, \underline{s}_{1}, \underline{s}_{2}, \overline{s}_{1}, \overline{s}_{2}, \widehat{l}_{1}, \widehat{l}_{2}, 0, \widehat{m}_{2}) .
        \end{equation*}

Since $(x, y_{1}, y_{2}) \in \mathbf{t}^{(- 1, j)} (\overline{s}_{1}, \underline{s}_{2}, \overline{s}_{1}, \overline{s}_{2})$, there exists also
        \begin{equation*}
            (0, \widetilde{l}_{2}, \widetilde{m}_{1}, \widetilde{m}_{2}) \in \mathbb{A}^{(- 1, j)} (x, y_{1}, y_{2}, \overline{s}_{1}, \underline{s}_{2}, \overline{s}_{1}, \overline{s}_{2})
        \end{equation*}
        such that
        \begin{equation*}
            \Phi (x, y_{1}, y_{2}, \overline{s}_{1}, \underline{s}_{2}, \overline{s}_{1}, \overline{s}_{2}) = \varphi (x, 0, y_{2}, \overline{s}_{1}, \underline{s}_{2}, \overline{s}_{1}, \overline{s}_{2}, \widetilde{l}_{1}, \widetilde{l}_{2}, \widetilde{m}_{1}, \widetilde{m}_{2}) .
        \end{equation*}
        Moreover $\mathbb{A}(x, y_{1}, y_{2}, \underline{s}_{1}, \underline{s}_{2}, \overline{s}_{1}, \overline{s}_{2})\subseteq
        \mathbb{A}(x, y_{1}, y_{2}, \overline{s}_{1}, \underline{s}_{2}, \overline{s}_{1}, \overline{s}_{2})$.
        For $\lambda := \frac{\widetilde{m}_{1}}{\widehat{l}_{1} + \widetilde{m}_{1}}$ consider the strategy
        \begin{equation*}
        (l_1,l_2,m_1,m_2):=\lambda (\widehat{l}_{1}, \widehat{l}_{2}, 0, \widehat{m}_{2})+(1-\lambda)(0, \widetilde{l}_{2}, \widetilde{m}_{1}, \widetilde{m}_{2})\in \mathbb{A}(x, y_{1}, y_{2}, \overline{s}_{1}, \underline{s}_{2}, \overline{s}_{1}, \overline{s}_{2}).
        \end{equation*}
        Then $\lambda \widehat{l}_{1}=(1-\lambda)\widetilde{m}_{1}$ and by Lemma \ref{Lem1.4} we have
        \begin{equation*}
        \Big( 0, \lambda \widehat{l}_{2} + (1 - \lambda ) \widetilde{l}_{2}, 0, \lambda \widehat{m}_{2} + (1 - \lambda ) \widetilde{m}_{2} \Big) \in \mathbb{A}(x, y_{1}, y_{2}, \overline{s}_{1}, \underline{s}_{2}, \overline{s}_{1}, \overline{s}_{2}) .
    \end{equation*}
    Such a strategy is also in $\mathbb{A}(x, y_{1}, y_{2}, \underline{s}_{1}, \underline{s}_{2}, \overline{s}_{1}, \overline{s}_{2})$. By Lemma \ref{concave} we have
    \begin{eqnarray*}
   && \varphi(x,  y_{1}, y_{2}, \underline{s}_{1}, \underline{s}_{2}, \overline{s}_{1}, \overline{s}_{2},  0, l_2, 0, m_2)=
  \varphi(x,  y_{1}, y_{2}, \underline{s}_{1}, \underline{s}_{2}, \overline{s}_{1}, \overline{s}_{2}, l_1, l_2, m_1, m_2) > \\
&& \lambda \varphi(x,  y_{1}, y_{2}, \underline{s}_{1}, \underline{s}_{2}, \overline{s}_{1}, \overline{s}_{2},\widehat{l}_{1}, \widehat{l}_{2}, 0, \widehat{m}_{2}) + (1-\lambda)\varphi(x,  y_{1}, y_{2}, \underline{s}_{1}, \underline{s}_{2}, \overline{s}_{1}, \overline{s}_{2}, 0, \widetilde{l}_{2}, \widetilde{m}_{1}, \widetilde{m}_{2})\geq \\
&&\Phi(x,  y_{1}, y_{2}, \underline{s}_{1}, \underline{s}_{2}, \overline{s}_{1}, \overline{s}_{2})
 \end{eqnarray*}
 and therefore the strategy $(\widehat{l}_{1}, \widehat{l}_{2}, 0, \widehat{m}_{2})$ is not optimal. We have a contradiction.
\end{proof}

\begin{proposition}\label{tobp7}
For $(\underline{s}_{1}, \underline{s}_{2}, \overline{s}_{1}, \overline{s}_{2}) \in \mathbb{D}$ and $i \in \{ - 1, 0, 1 \}$ we have
        \begin{equation}\label{bp7}
            \mathbf{t}^{(1, i)}(\underline{s}_{1}, \underline{s}_{2}, \overline{s}_{1}, \overline{s}_{2}) \subseteq \mathbf{t}^{(1, i)}(\overline{s}_{1}, \underline{s}_{2}, \overline{s}_{1}, \overline{s}_{2})
        \end{equation}
    \end{proposition}
 \begin{proof}
        Assume there exists $(x, y_{1}, y_{2}) \in \mathbf{t}^{(1, i)}(\underline{s}_{1}, \underline{s}_{2}, \overline{s}_{1}, \overline{s}_{2})\setminus \mathbf{t}^{(1, i)}(\overline{s}_{1}, \underline{s}_{2}, \overline{s}_{1}, \overline{s}_{2})$.
Then
\begin{eqnarray*}
&&(x, y_{1}, y_{2})  \in \bigcup_{(k_{1}, k_{2}) \in \{ - 1, 0, 1 \}^{2} \setminus \{ (1, i) \}} \mathbf{t}^{(k_{1}, k_{2})} (\overline{s}_{1}, \underline{s}_{2}, \overline{s}_{1}, \overline{s}_{2}) = \\
&& = \bigcup_{j \in \{ - 1, 0, 1 \}} \mathbf{t}^{(- 1, j)}(\overline{s}_{1}, \underline{s}_{2}, \overline{s}_{1}, \overline{s}_{2}) \cup \bigcup_{j \in \{ - 1, 0, 1 \}} \mathbf{t}^{(0, j)}(\overline{s}_{1}, \underline{s}_{2}, \overline{s}_{1}, \overline{s}_{2}) \\
&&\cup \bigcup_{j \in \{ - 1, 0, 1 \} \setminus \{ i \}} \mathbf{t}^{(1, j)}(\overline{s}_{1}, \underline{s}_{2}, \overline{s}_{1}, \overline{s}_{2}) .
\end{eqnarray*}
and by \eqref{bp4}
        \begin{equation*}
            (x, y_{1}, y_{2}) \notin \bigcup_{j \in \{ - 1, 0, 1 \}} \mathbf{t}^{(- 1, j)}(\overline{s}_{1}, \underline{s}_{2}, \overline{s}_{1}, \overline{s}_{2}) .
        \end{equation*}
By \eqref{bp3}
        \begin{equation*}
            (x, y_{1}, y_{2}) \notin \bigcup_{j \in \{ - 1, 0, 1 \}} \mathbf{t}^{(0, j)}(\overline{s}_{1}, \underline{s}_{2}, \overline{s}_{1}, \overline{s}_{2}),
        \end{equation*}
and finally by \eqref{bp6}
        \begin{equation*}
            (x, y_{1}, y_{2}) \notin \bigcup_{j \in \{ - 1, 0, 1 \} \setminus \{ i \}} \mathbf{t}^{(1, j)}(\overline{s}_{1}, \underline{s}_{2}, \overline{s}_{1}, \overline{s}_{2}),
        \end{equation*}
which together leads to a contradiction.
    \end{proof}
Combining \eqref{bp1} with Proposition \ref{tobp7}  we immediately obtain

    \begin{corollary}\label{tobp8}
    For $(\underline{s}_{1}, \underline{s}_{2}, \overline{s}_{1}, \overline{s}_{2}) \in \mathbb{D}$ and $i \in \{ - 1, 0, 1 \}$ we have
        \begin{equation}\label{bp8}
            \mathbf{t}^{(1, i)} (\underline{s}_{1}, \underline{s}_{2}, \overline{s}_{1}, \overline{s}_{2}) = \mathbf{t}^{(1, i)} (\overline{s}_{1}, \underline{s}_{2}, \overline{s}_{1}, \overline{s}_{2}) .
        \end{equation}
    \end{corollary}
Analogously one can obtain that

\begin{corollary}\label{tobp9}
    For $(\underline{s}_{1}, \underline{s}_{2}, \overline{s}_{1}, \overline{s}_{2}) \in \mathbb{D}$ and $i \in \{ - 1, 0, 1 \}$ we have
        \begin{equation}\label{bp9}
            \mathbf{t}^{(i, 1)} (\underline{s}_{1}, \underline{s}_{2}, \overline{s}_{1}, \overline{s}_{2}) = \mathbf{t}^{(i, 1)} (\underline{s}_{1}, \overline{s}_{2}, \overline{s}_{1}, \overline{s}_{2}) .
        \end{equation}
    \end{corollary}

 \section{Properties of no transactions zones}
 In this section we are going to show several properties of no transaction zone.
 \begin{proposition}\label{Prop7.1}
    For  $(\underline{s}_{1}, \underline{s}_{2}, \overline{s}_{1}, \overline{s}_{2}) \in \mathbb{D}$, ${s} \in [\underline{s}_{1}, \overline{s}_{1}]$ and $i \in \{ - 1, 0, 1 \}$ we have
    \begin{equation}\label{ino1}
       \mathbf{t}^{(0, i)} ({s}, \underline{s}_{2}, {s}, \overline{s}_{2}) \subseteq  \mathbf{t}^{(0, i)} (\underline{s}_{1}, \underline{s}_{2}, \overline{s}_{1}, \overline{s}_{2}) .
    \end{equation}
\end{proposition}
\begin{proof}
Assume that $(x,y_1,y_2)\in \mathbf{t}^{(0, i)} ({s}, \underline{s}_{2}, {s}, \overline{s}_{2})$. Then there exists $(l_1,l_2,m_1,m_2)\in \mathbb{A}^{(0, i)} (x, y_{1}, y_{2}, {s}, \underline{s}_{2}, {s}, \overline{s}_{2})$ such that
    \begin{equation*}
        \Phi (x, y_{1}, y_{2}, {s}, \underline{s}_{2}, {s}, \overline{s}_{2}) = \varphi (x, y_{1}, y_{2}, {s}, \underline{s}_{2}, {s}, \overline{s}_{2}, l_{1}, l_{2}, m_{1}, m_{2}) .
    \end{equation*}
    We also have that $(l_1,l_2,m_1,m_2)\in \mathbb{A} (x, y_{1}, y_{2}, \underline{s}_1, \underline{s}_{2}, \overline{s}_1, \overline{s}_{2})$. Therefore by Lemma \ref{usLem} we obtain
\begin{eqnarray*}
&&            \Phi (x, y_{1}, y_{2}, \underline{s}_{1}, \underline{s}_{2}, \overline{s}_{1}, \overline{s}_{2}) \leq \Phi (x, y_{1}, y_{2}, {s}, \underline{s}_{2}, {s}, \overline{s}_{2}) = \\
&& = \varphi (x, y_{1}, y_{2}, {s}, \underline{s}_{2}, {s}, \overline{s}_{2}, l_{1}, l_{2}, m_{1}, m_{2}) = \\
&& = \varphi (x, y_{1}, y_{2}, \underline{s}_{1}, \underline{s}_{2}, \overline{s}_{1}, \overline{s}_{2}, l_{1}, l_{2}, m_{1}, m_{2}) = \Phi (x, y_{1}, y_{2}, \underline{s}_{1}, \underline{s}_{2}, \overline{s}_{1}, \overline{s}_{2}).
\end{eqnarray*}
    which means that $(x,y_1,y_2)\in \mathbf{t}^{(0, i)} (x, y_{1}, y_{2}, \underline{s}_1, \underline{s}_{2}, \overline{s}_1, \overline{s}_{2})$.
\end{proof}
We now formulate our main result which is used in the construction of the shadow price.
\begin{theorem}\label{Thm7.1}
 For $(\underline{s}, \underline{s}_{2}, \overline{s}, \overline{s}_{2}) \in \mathbb{D}$ and
 $i \in \{ - 1, 0, 1 \}$ we have
 \begin{equation}\label{nodecom}
  \mathbf{t}^{(0, i)}(\underline{s}_{1}, \underline{s}_{2}, \overline{s}_{1}, \overline{s}_{2})
 = \bigcup_{s\in [\underline{s}_{1},\overline{s}_{1}]} \mathbf{t}^{(0, i)}({s}, \underline{s}_{2}, {s}, \overline{s}_{2}).
\end{equation}
\end{theorem}
\begin{proof} By Proposition \ref{Prop7.1} for $s\in [\underline{s}_{1},\overline{s}_{1}]$ we have \eqref{ino1}. Consequently
\begin{equation}
\bigcup_{s\in [\underline{s}_{1},\overline{s}_{1}]} \mathbf{t}^{(0, i)}({s}, \underline{s}_{2}, {s}, \overline{s}_{2})\subseteq \mathbf{t}^{(0, i)}(\underline{s}_{1}, \underline{s}_{2}, \overline{s}_{1}, \overline{s}_{2})
\end{equation}
and we have to show the inverse inclusion. Assume therefore that there is $(x,y_1,y_2)\in \mathbf{t}^{(0, i)}(\underline{s}_{1}, \underline{s}_{2}, \overline{s}_{1}, \overline{s}_{2})$ such that
\begin{equation}\label{notin}
(x,y_1,y_2)\notin \bigcup_{s\in [\underline{s}_{1},\overline{s}_{1}]}\mathbf{t}^{(0, i)}({s}, \underline{s}_{2}, {s}, \overline{s}_{2}).
\end{equation}\label{sit1}
Then for any $s\in [\underline{s}_{1},\overline{s}_{1}]$ there is $i(s)\neq i$ such that
\begin{equation}
(x,y_1,y_2) \in \mathbf{t}^{(0, i(s))}({s}, \underline{s}_{2}, {s}, \overline{s}_{2})
\end{equation}
or there is $j_1(s)\in\{-1,1\}$ and  $j_2(s)\in \{-1,0,1\}$ such that
\begin{equation}\label{sit2}
(x,y_1,y_2) \in \mathbf{t}^{(j_1(s), j_2(s))}({s}, \underline{s}_{2}, {s}, \overline{s}_{2}).
\end{equation}
In the case of \eqref{sit1} by \eqref{ino1} we have that $(x,y_1,y_2)\in \mathbf{t}^{(0, (s)}(\underline{s}_{1}, \underline{s}_{2}, \overline{s}_{1}, \overline{s}_{2})$ which is a contradiction to the fact that $i(s)\neq i$.

Assume therefore \eqref{sit2}. For $s=\overline{s}_1$ we have either $j_1(\overline{s}_1)=1$ or
$j_1(\overline{s}_1)=-1$. When $j_1(\overline{s}_1)=1$ by \eqref{bp1} we have
\begin{equation*}
\mathbf{t}^{(1, j_2(\overline{s}_1))} (\overline{s}_{1}, \underline{s}_{2}, \overline{s}_{1}, \overline{s}_{2}) \subseteq \mathbf{t}^{(1, j_2(\overline{s}_1))} (\underline{s}_{1}, \underline{s}_{2}, \overline{s}_{1}, \overline{s}_{2})
\end{equation*}
which contradicts the fact that  $(x,y_1,y_2)\in \mathbf{t}^{(0, i)}(\underline{s}_{1}, \underline{s}_{2}, \overline{s}_{1}, \overline{s}_{2})$. Consequently we have $j_1(\overline{s}_1)=-1$.
For $s=\underline{s}_1$ we again have either $j_1(\underline{s}_1)=1$ or
$j_1(\underline{s}_1)=-1$. When $j_1(\underline{s}_1)=-1$ by \eqref{prop1} we have
\begin{equation*}
\mathbf{t}^{(-1, j_2(\underline{s}_1))} (\underline{s}_{1}, \underline{s}_{2}, \underline{s}_{1}, \overline{s}_{2}) \subseteq \mathbf{t}^{(1, j_2(\underline{s}_1))} (\underline{s}_{1}, \underline{s}_{2}, \overline{s}_{1}, \overline{s}_{2})
\end{equation*}
which contradicts the fact that $(x,y_1,y_2)\in \mathbf{t}^{(0, i)}(\underline{s}_{1}, \underline{s}_{2}, \overline{s}_{1}, \overline{s}_{2})$. Therefore we have $j_1(\underline{s}_1)=1$. Summarizing we have that
\begin{equation}\label{waz}
(x,y_1,y_2)\in \mathbf{t}^{(-1, j_2(\overline{s}_1))} (\overline{s}_{1}, \underline{s}_{2}, \overline{s}_{1}, \overline{s}_{2}) \cap \mathbf{t}^{(1, j_2(\underline{s}_1))} (\underline{s}_{1}, \underline{s}_{2}, \underline{s}_{1}, \overline{s}_{2}).
\end{equation}
For fixed $(x,y_1,y_2,\underline{s}_2,\overline{s}_2)$ by Theorem \ref{select} there is a unique continuous function $ [\underline{s}_{1}, \overline{s}_{1}]\ni s \mapsto \left(\hat{l}_1(s),\hat{l}_2(s),\hat{m}_1(s),\hat{m}_2(s)\right)$ which forms an optimal investment strategy for the value function $\Phi(x,y_1,y_2,s,\underline{s}_2,s,\overline{s}_2)$. Notice that
by \eqref{waz} we have $\hat{l}_1(\overline{s}_1)=0$,  $\hat{m}_1(\overline{s}_1)>0$  and
$\hat{l}_1(\underline{s}_1)>0$,  $\hat{m}_1(\underline{s}_1)=0$. Since $\hat{l}_1(s)\hat{m}_1(s)=0$ we claim that there is $\hat{s}\in [\underline{s}_{1}, \overline{s}_{1}]$ such that $\hat{l}_1(\hat{s})=\hat{m}_1(\hat{s})=0$. In fact,
let $\hat{s}=\sup \left\{ s\in  [\underline{s}_{1}, \overline{s}_{1}]:\hat{l}_1(s)>0\right\}$.
By continuity of $s \mapsto \hat{l}_1(s)$   we clearly have that $\hat{l}_1(\hat{s})=0$. By continuity of $s \mapsto \hat{m}_1(s)$ we also have that  $\hat{m}_1(\hat{s})=0$. Consequently  $\hat{l}_1(\hat{s})=0$ and
$\hat{m}_1(\hat{s})=0$ and $(x,y_1,y_2)\in \mathbf{t}^{(0, j_2(\hat{s}))} (\hat{s}, \underline{s}_{2}, \hat{s}, \overline{s}_{2})$, but then by \eqref{ino1} we have
\begin{equation*}
(x,y_1,y_2)\in \mathbf{t}^{(0, j_2(\hat{s}))}((\underline{s}_{1}, \underline{s}_{2}, \overline{s}_{1}, \overline{s}_{2})
\end{equation*}
which may happen only when $j_2(\hat{s})=i$. But this means that
$(x,y_1,y_2) \in \mathbf{t}^{(0, i)}(\hat{s}, \underline{s}_{2}, \hat{s}, \overline{s}_{2})$
which is a contradiction to \eqref{notin}. This completes the proof.
\end{proof}
\begin{proposition}\label{Prop7.2}
     For $(\underline{s}, \underline{s}_{2}, \overline{s}, \overline{s}_{2}) \in \mathbb{D}$, $i, j \in \{ - 1, 0, 1 \}$ and $\widetilde{s} \in [\underline{s}, \overline{s}]$ we have

    \begin{equation}\label{7.1empty}
        \mathbf{t}^{(0, i)}(\overline{s}, \underline{s}_{2}, \overline{s}, \overline{s}_{2}) \cap \mathbf{t}^{(- 1, j)}(\widetilde{s}, \underline{s}_{2}, \widetilde{s}, \overline{s}_{2}) = \emptyset
    \end{equation}
and
\begin{equation}\label{7.2empty}
        \mathbf{t}^{(0, i)}(\underline{s}, \underline{s}_{2}, \underline{s}, \overline{s}_{2}) \cap \mathbf{t}^{(1, j)} (\widetilde{s}, \underline{s}_{2}, \widetilde{s}, \overline{s}_{2})=\emptyset.
    \end{equation}
    \end{proposition}

 \begin{proof}
 We prove first \eqref{7.1empty}. Assume that $ (x, y_{1}, y_{2}) \in \mathbf{t}^{(0, i)}(\overline{s}, \underline{s}_{2}, \overline{s}, \overline{s}_{2}) \cap \mathbf{t}^{(- 1, j)}(\widetilde{s}, \underline{s}_{2}, \widetilde{s}, \overline{s}_{2})$. Then there is $(\widetilde{l}_{1}, \widetilde{l}_{2}, \widetilde{m}_{1}, \widetilde{m}_{2}) \in \mathbb{A}^{(- 1, j)} (x, y_{1}, y_{2}, \widetilde{s}, \underline{s}_{2}, \widetilde{s}, \overline{s}_{2})$ such that $\widetilde{m}_{1}>0$ and
 \begin{equation}
        \Phi (x, y_{1}, y_{2}, \widetilde{s}, \underline{s}_{2}, \widetilde{s}, \overline{s}_{2}) = \varphi (x, y_{1}, y_{2}, \widetilde{s}, \underline{s}_{2}, \widetilde{s}, \overline{s}_{2},\widetilde{l}_{1}, \widetilde{l}_{2}, \widetilde{m}_{1}, \widetilde{m}_{2}) .
    \end{equation}
Since $\mathbb{A}^{(- 1, j)} (x, y_{1}, y_{2}, \widetilde{s}, \underline{s}_{2}, \widetilde{s}, \overline{s}_{2}) \subseteq \mathbb{A} (x, y_{1}, y_{2}, \overline{s}, \underline{s}_{2}, \overline{s}, \overline{s}_{2})$ we have that
\begin{equation}
\Phi (x, y_{1}, y_{2}, \widetilde{s}, \underline{s}_{2}, \widetilde{s}, \overline{s}_{2})\leq \Phi (x, y_{1}, y_{2}, \underline{s}, \underline{s}_{2}, \overline{s}, \overline{s}_{2}).
\end{equation}
Since also $\mathbb{A}^{(0, i)} (x, y_{1}, y_{2}, \overline{s}, \underline{s}_{2}, \overline{s}, \overline{s}_{2}) \subseteq \mathbb{A} (x, y_{1}, y_{2}, \widetilde{s}, \underline{s}_{2}, \widetilde{s}, \overline{s}_{2})$ we have (taking into account that $(x, y_{1}, y_{2}) \in \mathbf{t}^{(0, i)}(\overline{s}, \underline{s}_{2}, \overline{s}, \overline{s}_{2})$) that
\begin{equation}\label{eqqw}
\Phi (x, y_{1}, y_{2}, \underline{s}, \underline{s}_{2}, \overline{s}, \overline{s}_{2})\leq \Phi (x, y_{1}, y_{2}, \widetilde{s}, \underline{s}_{2}, \widetilde{s}, \overline{s}_{2})
\end{equation}
and consequently  $\Phi (x, y_{1}, y_{2}, \widetilde{s}, \underline{s}_{2}, \widetilde{s}, \overline{s}_{2})=\Phi (x, y_{1}, y_{2}, \underline{s}, \underline{s}_{2}, \overline{s}, \overline{s}_{2})$. Moreover there is a strategy  $ (l_{1}, l_{2}, m_{1}, m_{2}) \in \mathbb{A}^{(0, i)} (x, y_{1}, y_{2}, \overline{s}, \underline{s}_{2}, \overline{s}, \overline{s}_{2})$ with $l_1=m_1=0$ such that
\begin{eqnarray}
&&\Phi (x, y_{1}, y_{2}, \underline{s}, \underline{s}_{2}, \overline{s}, \overline{s}_{2})=\varphi (x, y_{1}, y_{2}, \underline{s}, \underline{s}_{2}, \overline{s}, \overline{s}_{2},{l}_{1}, {l}_{2}, {m}_{1}, {m}_{2})=\nonumber \\
&&\varphi (x, y_{1}, y_{2}, \widetilde{s}, \underline{s}_{2}, \widetilde{s}, \overline{s}_{2},{l}_{1}, {l}_{2}, {m}_{1}, {m}_{2})
\end{eqnarray}
and we have two optimal strategies for  $\Phi (x, y_{1}, y_{2}, \widetilde{s}, \underline{s}_{2}, \widetilde{s}, \overline{s}_{2})$ which contradicts Proposition \ref{Prop2.3}.

\noindent
The proof of \eqref{7.2empty} is similar. Assume that $ (x, y_{1}, y_{2}) \in \mathbf{t}^{(0, i)}(\overline{s}, \underline{s}_{2}, \overline{s}, \overline{s}_{2})\cap \mathbf{t}^{(- 1, j)}(\widetilde{s}, \underline{s}_{2}, \widetilde{s}, \overline{s}_{2})$. Then there is $(\widetilde{l}_{1}, \widetilde{l}_{2}, \widetilde{m}_{1}, \widetilde{m}_{2}) \in \mathbb{A}^{(1, j)} (x, y_{1}, y_{2}, \widetilde{s}, \underline{s}_{2}, \widetilde{s}, \overline{s}_{2})$ such that $\widetilde{l}_{1}>0$ and
 \begin{equation}
        \Phi (x, y_{1}, y_{2}, \widetilde{s}, \underline{s}_{2}, \widetilde{s}, \overline{s}_{2}) = \varphi (x, y_{1}, y_{2}, \widetilde{s}, \underline{s}_{2}, \widetilde{s}, \overline{s}_{2},\widetilde{l}_{1}, \widetilde{l}_{2}, \widetilde{m}_{1}, \widetilde{m}_{2}) .
    \end{equation}
 Since $\mathbb{A}^{(1, j)} (x, y_{1}, y_{2}, \widetilde{s}, \underline{s}_{2}, \widetilde{s}, \overline{s}_{2}) \subseteq \mathbb{A} (x, y_{1}, y_{2}, \overline{s}, \underline{s}_{2}, \overline{s}, \overline{s}_{2})$ we have that
\begin{equation}
\Phi (x, y_{1}, y_{2}, \widetilde{s}, \underline{s}_{2}, \widetilde{s}, \overline{s}_{2})\leq \Phi (x, y_{1}, y_{2}, \underline{s}, \underline{s}_{2}, \overline{s}, \overline{s}_{2}).
\end{equation}
Now by \eqref{eqqw} we obtain that $\Phi (x, y_{1}, y_{2}, \widetilde{s}, \underline{s}_{2}, \widetilde{s}, \overline{s}_{2})= \Phi (x, y_{1}, y_{2}, \underline{s}, \underline{s}_{2}, \overline{s}, \overline{s}_{2})$ and there is $ (l_{1}, l_{2}, m_{1}, m_{2}) \in \mathbb{A}^{(0, i)} (x, y_{1}, y_{2}, \overline{s}, \underline{s}_{2}, \overline{s}, \overline{s}_{2})$ with $l_1=m_1=0$ such that
\begin{eqnarray}
&&\Phi (x, y_{1}, y_{2}, \underline{s}, \underline{s}_{2}, \overline{s}, \overline{s}_{2})=\varphi (x, y_{1}, y_{2}, \underline{s}, \underline{s}_{2}, \overline{s}, \overline{s}_{2},{l}_{1}, {l}_{2}, {m}_{1}, {m}_{2})=\nonumber \\
&&\varphi (x, y_{1}, y_{2}, \widetilde{s}, \underline{s}_{2}, \widetilde{s}, \overline{s}_{2},{l}_{1}, {l}_{2}, {m}_{1}, {m}_{2}),
\end{eqnarray}
which means that we have two optimal strategies for $\Phi (x, y_{1}, y_{2}, \widetilde{s}, \underline{s}_{2}, \widetilde{s}, \overline{s}_{2})$ and we have a contradiction to Proposition \ref{Prop2.3}.
\end{proof}

\begin{proposition}\label{Prop9.3}
For $(\underline{s}, \underline{s}_{2}, \overline{s}, \overline{s}_{2}) \in \mathbb{D}$ and $\widetilde{s} \in [\underline{s}, \overline{s}]$,  $i, j \in \{ - 1, 0, 1 \}$  such that $i \neq j$ we have
    \begin{equation}\label{P9.31}
        \mathbf{t}^{(0, i)}(\underline{s}, \underline{s}_{2}, \underline{s}, \overline{s}_{2}) \cap \mathbf{t}^{(0, j)}(\widetilde{s}, \underline{s}_{2}, \widetilde{s}, \overline{s}_{2}) = \emptyset.
    \end{equation}
\begin{equation}\label{P9.32}
        \mathbf{t}^{(0, i)}(\overline{s}, \underline{s}_{2}, \overline{s}, \overline{s}_{2}) \cap \mathbf{t}^{(0, j)}(\widetilde{s}, \underline{s}_{2}, \widetilde{s}, \overline{s}_{2}) = \emptyset.
    \end{equation}
    \end{proposition}
    \begin{proof} We prove only \eqref{P9.31} since the proof of \eqref{P9.32} is almost the same. Assume that
     $(x, y_{1}, y_{2}) \in \mathbf{t}^{(0, i)}(\underline{s}, \underline{s}_{2}, \underline{s}, \overline{s}_{2}) \cap \mathbf{t}^{(0, j)}(\widetilde{s}, \underline{s}_{2}, \widetilde{s}, \overline{s}_{2})$.

     Since $ \mathbb{A}^{(0, k)}(x, y_{1}, y_{2}, \underline{s}, \underline{s}_{2}, \underline{s}, \overline{s}_{2}) = \mathbb{A}^{(0, k)}(x, y_{1}, y_{2}, \widetilde{s}, \underline{s}_{2}, \widetilde{s}, \overline{s}_{2})$ for $k=i,j$
     we have that
     $\Phi (x, y_{1}, y_{2}, \underline{s}, \underline{s}_{2}, \underline{s}, \overline{s}_{2}) = \Phi (x, y_{1}, y_{2}, \widetilde{s}, \underline{s}_{2}, \widetilde{s}, \overline{s}_{2})$. Consequently there are:

      $(l_{1}, l_{2}, m_{1}, m_{2}) \in \mathbb{A}^{(0, i)} (x, y_{1}, y_{2}, \underline{s}, \underline{s}_{2}, \underline{s}, \overline{s}_{2})$ and
     $(\widetilde{l}_{1}, \widetilde{l}_{2}, \widetilde{m}_{1}, \widetilde{m}_{2}) \in \mathbb{A}^{(0, j)}(x, y_{1}, y_{2}, \widetilde{s}, \underline{s}_{2}, \widetilde{s}, \overline{s}_{2})$
     such that
     \begin{eqnarray}
     &&\Phi (x, y_{1}, y_{2}, \underline{s}, \underline{s}_{2}, \underline{s}, \overline{s}_{2})= \varphi (x, y_{1}, y_{2}, \underline{s}, \underline{s}_{2}, \underline{s}, \overline{s}_{2},{l}_{1}, {l}_{2}, {m}_{1},{m}_{2})= \nonumber \\
     &&  \varphi (x, y_{1}, y_{2}, \widetilde{s}, \underline{s}_{2}, \widetilde{s}, \overline{s}_{2},\widetilde{l}_{1}, \widetilde{l}_{2}, \widetilde{m}_{1}, \widetilde{m}_{2})= \varphi (x, y_{1}, y_{2}, \underline{s}, \underline{s}_{2}, \underline{s}, \overline{s}_{2},\widetilde{l}_{1}, \widetilde{l}_{2}, \widetilde{m}_{1}, \widetilde{m}_{2}),
     \end{eqnarray}
     and we have two different optimal strategies for $\Phi (x, y_{1}, y_{2}, \underline{s}, \underline{s}_{2}, \underline{s}, \overline{s}_{2})$ a contradiction to Proposition \ref{Prop2.3}.
     \end{proof}

\begin{proposition}\label{Prop9.5}
  Assume $(\underline{s}_{1}, \underline{s}, \overline{s}_{1}, \overline{s}) \in \mathbb{D}$, $i \in \{ - 1, 0, 1 \}$ and $\widetilde{s} \in [\underline{s}, \overline{s}]$. Then
  \begin{equation}\label{P9.11}
        \mathbf{t}^{(0, i)}(\underline{s}, \underline{s}_{2}, \underline{s}, \overline{s}_{2}) \cap \mathbf{t}^{(0, i)}(\overline{s}, \underline{s}_{2}, \overline{s}, \overline{s}_{2}) \subseteq \mathbf{t}^{(0, i)}(\widetilde{s}, \underline{s}_{2}, \widetilde{s}, \overline{s}_{2}) .
    \end{equation}
    and
    \begin{equation}\label{P9.12}
        \mathbf{t}^{(i, 0)}(\underline{s}_{1}, \underline{s}, \overline{s}_{1}, \underline{s}) \cap \mathbf{t}^{(i, 0)}(\underline{s}_{1}, \overline{s}, \overline{s}_{1}, \overline{s}) \subseteq \mathbf{t}^{(i, 0)}(\underline{s}_{1}, \widetilde{s}, \overline{s}_{1}, \widetilde{s}) .
    \end{equation}
\end{proposition}
\begin{proof}
Notice that
\begin{eqnarray*}
           && (x, y_{1}, y_{2}) \in \big( \mathbb{R}_{+} \times \mathbb{R}_{+}^{2} \big) \setminus \mathbf{t}^{(0, i)}(\widetilde{s}, \underline{s}_{2}, \widetilde{s}, \overline{s}_{2}) =\\
            && = \bigcup_{j \in \{ - 1, 0, 1 \}} \mathbf{t}^{(- 1, j)}(\widetilde{s}, \underline{s}_{2}, \widetilde{s}, \overline{s}_{2}) \cup \bigcup_{j \in \{ - 1, 0, 1 \} \setminus \{ i \}} \mathbf{t}^{(0, j)}(\widetilde{s}, \underline{s}_{2}, \widetilde{s}, \overline{s}_{2}) \cup \bigcup_{j \in \{ - 1, 0, 1 \}} \mathbf{t}^{(1, j)}(\widetilde{s}, \underline{s}_{2}, \widetilde{s}, \overline{s}_{2}) .
    \end{eqnarray*}
    Therefore by Propositions \ref{Prop7.1} and \ref{Prop7.2} we immediately have \eqref{P9.11}. The proof of \eqref{P9.12} is by analogy to that of \eqref{P9.11}.
    \end{proof}

\section{Construction of the shadow price}

Let $(\underline{s}_{1}, \underline{s}_{2}, \overline{s}_{1}, \overline{s}_{2}) \in \mathbb{D}$.
From Lemma \ref{disjoint} we know that for every $(x, y_{1}, y_{2}) \in \mathbb{R}_{+} \times \mathbb{R}_{+}^{2}$ there exists a unique pair $(i_{1}, i_{2}) \in \{ - 1, 0, 1 \}^{2}$ such that $(x, y_{1}, y_{2}) \in \mathbf{t}^{(i_{1}, i_{2})}(\underline{s}_{1}, \underline{s}_{2}, \overline{s}_{1}, \overline{s}_{2})$.
Therefore for  $(x, y_{1}, y_{2}) \in \mathbb{R}_{+} \times \mathbb{R}_{+}^{2}$ we can define
\begin{equation*}
    \mathrm{I}_{1}(x, y_{1}, y_{2}), \mathrm{I}_{2}(x, y_{1}, y_{2}) \in \{ - 1, 0, 1 \}
\end{equation*}
such that
\begin{equation}
    (x, y_{1}, y_{2}) \in \mathbf{t}^{\big( \mathrm{I}_{1}(x, y_{1}, y_{2}), \mathrm{I}_{2}(x, y_{1}, y_{2}) \big)}(\underline{s}_{1}, \underline{s}_{2}, \overline{s}_{1}, \overline{s}_{2}) .
\end{equation}

From Theorem \ref{Thm7.1}  for every $i \in \{ - 1, 0, 1 \}$ and for every $(x, y_{1}, y_{2}) \in \mathbf{t}^{(0, i)}(\underline{s}_{1}, \underline{s}_{2}, \overline{s}_{1}, \overline{s}_{2})$ there exists an $\widetilde{s} \in [\underline{s}_{1}, \overline{s}_{1}]$ such that $(x, y_{1}, y_{2}) \in \mathbf{t}^{(0,i)}(\widetilde{s}, \underline{s}_{2}, \widetilde{s}, \overline{s}_{2})$.

Therefore, for every $(x, y_{1}, y_{2}) \in \mathbb{R}_{+} \times \mathbb{R}_{+}^{2}$ we can define the set
\begin{equation}\label{sha1}
    {\cal  S}_{1}(x, y_{1}, y_{2}) :=             \begin{cases}
                \{ \underline{s}_{1} \} \quad & \mbox{if $\mathrm{I}_{1}(x, y_{1}, y_{2}) = - 1$} \\
                \Big\{ s \in [\underline{s}_{1}, \overline{s}_{1}] : & \mbox{when $(x, y_{1}, y_{2}) \in \mathbf{t}^{\big( 0, \mathrm{I}_{2}(x, y_{1}, y_{2}) \big)} (s, \underline{s}_{2}, s, \overline{s}_{2})$} \\
                 &\mbox{and $\mathrm{I}_{1}(x, y_{1}, y_{2}) = 0\Big\}$} \\
                \{ \overline{s}_{1} \} \quad & \mbox{if $\mathrm{I}_{1}(x, y_{1}, y_{2}) = 1$}
            \end{cases}.
\end{equation}
and this is a nonempty set. Furthermore we can show that

\begin{proposition}\label{path-connected}
    For every $(x, y_{1}, y_{2}) \in \mathbb{R}_{+} \times \mathbb{R}_{+}^{2}$ the set ${\cal S}_{1}(x, y_{1}, y_{2})$ is path-connected.
\end{proposition}

\begin{proof}
    Let $(x, y_{1}, y_{2}) \in \mathbb{R}_{+} \times \mathbb{R}_{+}^{2}$ be arbitrary.
It is clear that we have to consider only  the case, when $\mathrm{I}_{1}(x, y_{1}, y_{2}) = 0$

    Let $s_{1}, s_{2} \in {\it S}_{1}(x, y_{1}, y_{2})$ be arbitrary. Without loss of generality we will assume that $s_{1} \leq s_{2}$ and show that for every $s \in [s_{1}, s_{2}]$ is in  $ {\cal S}_{1}(x, y_{1}, y_{2})$.
Since  $s_{1}, s_{2} \in {\cal S}_{1}(x, y_{1}, y_{2})$ we have
    \begin{equation*}
        (x, y_{1}, y_{2}) \in \mathbf{t}^{\big( 0, \mathrm{I}_{2}(x, y_{1}, y_{2}) \big)} (s_{1}, \underline{s}_{2}, s_{1}, \overline{s}_{2}) \cap \mathbf{t}^{\big( 0, \mathrm{I}_{2}(x, y_{1}, y_{2}) \big)} (s_{2}, \underline{s}_{2}, s_{2}, \overline{s}_{2}) .
    \end{equation*}
    By \eqref{P9.11} for every $s \in [\underline{s}_{1}, \overline{s}_{2}]$ we have that
    $(x, y_{1}, y_{2}) \in \mathbf{t}^{\big( 0, \mathrm{I}_{2}(x, y_{1}, y_{2}) \big)} (s, \underline{s}_{2}, s, \overline{s}_{2}),$
which means that $s \in {\cal S}_{1}(x, y_{1}, y_{2})$ for every $s \in [s_{1}, s_{2}]$. This completes the proof.
  \end{proof}
Directly from the above Proposition we obtain
\begin{corollary}
    For every $(x, y_{1}, y_{2}) \in \mathbb{R}_{+} \times \mathbb{R}_{+}^{2}$ the set ${\cal S}_{1}(x, y_{1}, y_{2})$ is either an interval or a singleton.
\end{corollary}

We study further properties of the set ${\cal S}$.

\begin{proposition}\label{closed}
    For every $(x, y_{1}, y_{2}) \in \mathbb{R}_{+} \times \mathbb{R}_{+}^{2}$ the set ${\cal S}_{1}(x, y_{1}, y_{2})$ is closed.
\end{proposition}

\begin{proof}
Since when $\mathrm{I}_{1}(x, y_{1}, y_{2}) \in \{ - 1, 1 \}$ the set ${\cal S}_{1}(x, y_{1}, y_{2})$ is a singleton we have to consider only the case  $\mathrm{I}_{1}(x, y_{1}, y_{2}) = 0$. Assume sequence $s^n\in {\cal S}_{1}(x, y_{1}, y_{2})$ converges
to some $\widetilde{s} \in [\underline{s}_{1}, \overline{s}_{1}]$.
    Let $\widehat{u} : {\it S}_{1}(x, y_{1}, y_{2}) \longrightarrow \mathbb{R}_{+}^{2} \times \mathbb{R}_{+}^{2}$ be such a function that for every $s \in {\cal S}_{1}(x, y_{1}, y_{2})$ we have $\widehat{u}(s) \in \overleftrightarrow{\mathbb{A}}(x, y_{1}, y_{2}, s, \underline{s}_{2}, s, \overline{s}_{2})$ and
    $\Phi (x, y_{1}, y_{2}, s, \underline{s}_{2}, s, \overline{s}_{2}) = \varphi \big(x, y_{1}, y_{2}, s, \underline{s}_{2}, s, \overline{s}_{2}, \widehat{u}(s) \big).$
    From Corollary $\ref{Corollary on the continuity of the function widehat}$ we know that the function $\widehat{u}$ is continuous and consequently $\widehat{u}(s^{n}) \longrightarrow \widehat{u}(\widetilde{s})$ as $n \longrightarrow \infty$.
Since $\mathrm{I}_{1}(x, y_{1}, y_{2}) = 0$, then for every $n \in \mathbb{N}$ we have that $\widehat{u}_{1}(s^{n}) = 0$ and $\widehat{u}_{3}(s^{n}) = 0$. Therefore also $\widehat{u}_{1}(\widetilde{s}) = 0$ and $\widehat{u}_{3}(\widetilde{s}) = 0$ and
    \begin{equation*}
        (x, y_{1}, y_{2}) \in \bigcup_{j \in \{ - 1, 0, 1 \}} \mathbf{t}^{(0, j)}(\widetilde{s}, \underline{s}_{2}, \widetilde{s}, \overline{s}_{2}),
    \end{equation*}
    which means that $ (x, y_{1}, y_{2}) \in  \mathbf{t}^{(0, j)}(\widetilde{s}, \underline{s}_{2}, \widetilde{s}, \overline{s}_{2})$ for some $j \in \{ - 1, 0, 1 \}$.
From Theorem \ref{Thm7.1} we have that
$(x, y_{1}, y_{2}) \in \mathbf{t}^{\big( 0, j \big)} (\underline{s}_1, \underline{s}_{2}, \overline{s}_1, \overline{s}_{2}),$ and
therefore $j=\mathrm{I}_{2}(x, y_{1}, y_{2})$. Consequently  $\widetilde{s} \in {\it S}_{1}(x, y_{1}, y_{2})$.
\end{proof}

For every $(x, y_{1}, y_{2}) \in \mathbb{R}_{+} \times \mathbb{R}_{+}^{2}$ define the function $\widetilde{S}_{1}: \mathbb{R}_{+} \times \mathbb{R}_{+}^{2} \longrightarrow [\underline{s}_{1}, \overline{s}_{1}]$ by the formula
\begin{equation}\label{shad1}
    \widetilde{S}_{1}(x, y_{1}, y_{2}) := \inf {\cal S}_{1}(x, y_{1}, y_{2}).
\end{equation}

From Proposition \ref{closed} it is clear that
$\widetilde{S}_{1}$ takes values in the set ${\cal S}_{1}(x, y_{1}, y_{2})$. In what follows we want to show measurability of
$ \widetilde{S}_{1}(x, y_{1}, y_{2})$. We have

\begin{proposition}\label{fmeasurable}
      Assume that function $f_{k} : \mathbb{R}_{+} \times \mathbb{R}_{+}^{2} \longleftrightarrow \{ 0 \} \cup [\underline{s}_{1}, \overline{s}_{1}]$ for $k \in \{ - 1, 0, 1 \}$ is given  by the formula
    \begin{equation}
        f_{k}(x, y_{1}, y_{2}) = \widetilde{S}_{1}(x, y_{1}, y_{2}) \cdot \mathrm{1}_{\big\{ \mathrm{I}_{1}(x, y_{1}, y_{2}) = 0 \big\}}(x, y_{1}, y_{2}) \cdot \mathrm{1}_{\big\{ \mathrm{I}_{2}(x, y_{1}, y_{2}) = k \big\}}(x, y_{1}, y_{2}) .
    \end{equation}
    Then $f_{k}$ is a $\mathcal{B}(\mathbb{R}_{+} \times \mathbb{R}_{+}^{})$-measurable function.
\end{proposition}

\begin{proof}
    For $a > \overline{s}_{1}$ we have that $f_{k}^{- 1}\big[ [- \infty , a] \big] = \mathbb{R}_{+} \times \mathbb{R}_{+}^{2} \in \mathcal{B}(\mathbb{R}_{+} \times \mathbb{R}_{+}^{2})$ .

      Consider now $a \in [\underline{s}_{1}, \overline{s}_{1}]$. From Theorem \ref{Thm7.1} and then Lemma \ref{disjoint} we have
    \begin{eqnarray*}
    &&f_{k}^{- 1} \big[[\underline{s}_{1}, a] \big] = \big\{ (x, y_{1}, y_{2}) \in \mathbb{R}_{+} \times \mathbb{R}_{+}^{2} : \quad \underline{s}_{1} \leq f_{k}(x, y_{1}, y_{2}) \leq a \big\} = \\
            && = \Big\{ (x, y_{1}, y_{2}) \in \mathbb{R}_{+} \times \mathbb{R}_{+}^{2} : \quad \underline{s}_{1} \leq \widetilde{S}_{1}(x, y_{1}, y_{2}) \cdot \mathrm{1}_{\big\{ \mathrm{I}_{1}(x, y_{1}, y_{2}) = 0 \big\}}(x, y_{1}, y_{2})\\
            && \cdot \mathrm{1}_{\big\{ \mathrm{I}_{2}(x, y_{1}, y_{2}) = k \big\}}(x, y_{1}, y_{2}) \leq a \Big\} = \\
            && = \Big\{ (x, y_{1}, y_{2}) \in \mathbb{R}_{+} \times \mathbb{R}_{+}^{2} : \quad \underline{s}_{1} \leq \widetilde{S}_{1}(x, y_{1}, y_{2}) \leq a \Big\} \cap \mathbf{t}^{(0, k)}(\underline{s}_{1}, \underline{s}_{2}, \overline{s}_{1}, \overline{s}_{2}) = \\
            && = \Big\{ (x, y_{1}, y_{2}) \in \mathbb{R}_{+} \times \mathbb{R}_{+}^{2} : \quad \underset{\widetilde{s} \in [\underline{s_{1}}, a]}{\exists} \, (x, y_{1}, y_{2}) \in \mathbf{t}^{(0, k)}(\widetilde{s}, \underline{s}_{2}, \widetilde{s}, \overline{s}_{2}) \Big\} \cap \mathbf{t}^{(0, k)}(\underline{s}_{1}, \underline{s}_{2}, \overline{s}_{1}, \overline{s}_{2}) = \\
            &&=\Bigg[ \bigcup_{\widetilde{s} \in [\underline{s}_{1}, a]} \mathbf{t}^{(0, k)}(\widetilde{s}, \underline{s}_{2}, \widetilde{s}, \overline{s}_{2}) \Bigg] \cap \Bigg[ \bigcup_{\widetilde{s} \in [\underline{s}_{1}, \overline{s}_{1}]} \mathbf{t}^{(0, k)}(\widetilde{s}, \underline{s}_{2}, \widetilde{s}, \overline{s}_{2}) \Bigg] = \\
            &&=  \bigcup_{\widetilde{s} \in [\underline{s}_{1}, a]} \mathbf{t}^{(0, k)}(\widetilde{s}, \underline{s}_{2}, \widetilde{s}, \overline{s}_{2}) = \mathbf{t}^{(0, k)}(\underline{s}_{1}, \underline{s}_{2}, a, \overline{s}_{2}) \in \mathcal{B}(\mathbb{R}_{+} \times \mathbb{R}_{+}^{2}) .
    \end{eqnarray*}

    For $a \in [0, \underline{s}_{1})$ using again Lemma \ref{disjoint} we obtain
    \begin{equation*}
        f_{k}^{- 1} \big[ [- \infty , a] \big] = f_{k}^{- 1} \big[ \{ 0 \} \big] = \bigcup_{(i, j) \in \{ - 1, 0, 1 \}^{2} \setminus \{ (0, k) \}} \mathbf{t}^{(i,j)}(\underline{s}_{1}, \underline{s}_{2}, \overline{s}_{1}, \overline{s}_{2}) \in \mathcal{B}(\mathbb{R}_{+} \times \mathbb{R}_{+}^{2}) .
    \end{equation*}
Finally for  $a < 0$ we have $ f_{k}^{- 1} \big[ [- \infty, a] \big] = \emptyset \in \mathcal{B}(\mathbb{R}_{+} \times \mathbb{R}_{+}^{2}).$ Therefore for $a \in ( -\infty , \overline{s}_{1}]$ we obtain that
    $f_{k}^{- 1} \big[ [ -\infty , a] \big] \in \mathcal{B}(\mathbb{R}_{+} \times \mathbb{R}_{+}^{2})$.
Summarizing, $f_{k}$ is $\mathcal{B}(\mathbb{R}_{+} \times \mathbb{R}_{+}^{2})$-measurable.
\end{proof}

We now complete the construction of shadow price $\widetilde{S}_{1}(x, y_{1}, y_{2})$.
\begin{corollary}\label{shameas}
    Function $\widetilde{S}_{1}$ is $\mathcal{B}(\mathbb{R}_{+} \times \mathbb{R}_{+}^{2})$-measurable.
\end{corollary}

\begin{proof}
For very $(x, y_{1}, y_{2}) \in \mathbb{R}_{+} \times \mathbb{R}_{+}^{2}$ we have that
    \begin{eqnarray*}
           && \widetilde{S}_{1}(x, y_{1}, y_{2})  = \underline{s}_{1} \cdot \mathrm{1}_{\big\{ \mathrm{I}_{1}(x, y_{1}, y_{2}) = - 1 \big\}}(x, y_{1}, y_{2}) + \overline{s}_{1} \cdot \mathrm{1}_{\big\{ \mathrm{I}_{1}(x, y_{1}, y_{2}) = 1 \big\}}(x, y_{1}, y_{2}) + \\
            && \ \ \ \ \ \ + \sum_{k = - 1}^{1} \widetilde{S}_{1}(x, y_{1}, y_{2}) \cdot \mathrm{1}_{\big\{ \mathrm{I}_{1}(x, y_{1}, y_{2}) = 0 \big\}}(x, y_{1}, y_{2}) \cdot \mathrm{1}_{\big\{ \mathrm{I}_{2}(x, y_{1}, y_{2}) = k \big\}}(x, y_{1}, y_{2}) .
    \end{eqnarray*}
    Thus, from Proposition \ref{fmeasurable}  function $\widetilde{S}_{1}$ is $\mathcal{B}(\mathbb{R}_{+} \times \mathbb{R}_{+}^{2})$-measurable as the sum of $\mathcal{B}(\mathbb{R}_{+} \times \mathbb{R}_{+}^{2})$- measurable functions.
\end{proof}
By analogy  for every $(x, y_{1}, y_{2}) \in \mathbf{t}^{\big( \mathrm{I}_{1}(x, y_{1}, y_{2}), \mathrm{I}_{2}(x, y_{1}, y_{2}) \big)}(\underline{s}_{1}, \underline{s}_{2}, \overline{s}_{1}, \overline{s}_{2})$  we can define the set
\begin{equation}
    {\cal  S}_{2}(x, y_{1}, y_{2}) :=             \begin{cases}
                \{ \underline{s}_{2} \} \quad & \mbox{if $\mathrm{I}_{2}(x, y_{1}, y_{2}) = - 1$} \\
                \Big\{ s \in [\underline{s}_{2}, \overline{s}_{2}] : & \mbox{when $(x, y_{1}, y_{2}) \in \mathbf{t}^{\big(\mathrm{I}_{1}(x, y_{1}, y_{2}),0 \big)} (\underline{s}_{1}, s, \overline{s}_{1},s)$} \\
                 &\mbox{and $\mathrm{I}_{2}(x, y_{1}, y_{2}) = 0\Big\}$} \\
                \{ \overline{s}_{2} \} \quad & \mbox{if $\mathrm{I}_{2}(x, y_{1}, y_{2}) = 1$}
            \end{cases}.
\end{equation}
and following \eqref{shad1} define
\begin{equation}\label{shad2}
    \widetilde{S}_{2}(x, y_{1}, y_{2}) := \inf {\cal S}_{2}(x, y_{1}, y_{2}).
\end{equation}
We have
\begin{theorem}\label{sshadow}
The pair $\left(\widetilde{S}_{1}(x, y_{1}, y_{2}), \widetilde{S}_{2}(x, y_{1}, y_{2})\right)$ defined in \eqref{shad1} and \eqref{shad2} form shadow price for the static problem i.e. we have with $\Phi$ defined in \eqref{Phi}
\begin{eqnarray}\label{shadd}
&&\Phi (x, y_{1}, y_{2}, \underline{s}_{1}, \underline{s}_{2}, \overline{s}_{1}, \overline{s}_{2}) := \nonumber \\
&&\sup_{(l_{1}, l_{2}, m_{1}, m_{2}) \in \mathbb{A}(x, y_{1}, y_{2}, \widetilde{S}_{1}, \widetilde{S}_{2}, \widetilde{S}_{1}, \widetilde{S}_{2})} \varphi (x, y_{1}, y_{2}, \widetilde{S}_{1}, \widetilde{S}_{2}, \widetilde{S}_{1}, \widetilde{S}_{2}, l_{1}, l_{2}, m_{1}, m_{2}) .
\end{eqnarray}
where to simplify notations we write $\widetilde{S}_i$ for $\widetilde{S}_{2}(x, y_{1}, y_{2})$.
\end{theorem}
\begin{proof} We follow the arguments of Propositions \ref{path-connected}, \ref{closed}, \ref{fmeasurable} and Corollary \ref{shameas}.
\end{proof}

\section{Dynamic shadow price for two asset case}
In this section we consider two assets with bid  $(\underline{S}_1(t),\underline{S}_2(t))$ and ask $(\overline{S}_1(t),\overline{S}_2(t))$ prices over time discrete time interval $[0,T]$. Consider regular conditional laws $\mathbb{P}_t^{i,j}$ (which exist by Theorem 6.3 of \cite{Kal}) defined for $i,j\in\left\{-1,1\right\}$ as
\begin{equation}\label{Ass1}
\mathbb{P}_t^{i,j}(\cdot)=\mathbb{P}\left((S_1^i(t+1),S_2^j(t+1))\in \cdot|{\cal F}_t\right),
\end{equation}
with the notation $S_k^i(t+1)=\underline{S}_{k}(t+1)$ when $i=-1$ and $S_k^i(t+1)=\overline{S}_{k}(t+1)$ when $i=1$ with $k=1,2$.
We shall assume that
\begin{assumption}\label{ass1}
Conditional laws $\mathbb{P}_t^{i,j}$  for $i,j\in\left\{-1,1\right\}$ are non degenerate, which means that
\begin{equation}
\mathbb{P}\left(a_1 S_1^i(t+1)+a_2 S_2^j(t+1) = a_0 |{\cal F}_t\right)=0
\end{equation}
for any ${\cal F}_t$ measurable random variables $a_0,a_1,a_2$ and $t=0,1,\ldots,T-1$. 
\end{assumption}
Define now using regular conditional probabilities  the following system of random Bellman equations for $(x,y_1,y_2,\underline{s}_{1}, \underline{s}_{2}, \overline{s}_{1}, \overline{s}_{2}) \in \mathbb{R}_{+} \times \mathbb{R}_{+}^{2} \times \mathbb{D}$

\begin{equation}\label{bel1}
w_T(x,y_1,y_2,\underline{s}_{1}, \underline{s}_{2}, \overline{s}_{1}, \overline{s}_{2}):=U(x+y_1 \underline{s}_{1} + y_2 \underline{s}_{2}),
\end{equation}
\begin{eqnarray}\label{bel2}
&& w_{T-1}(x,y_1,y_2,\underline{s}_{1}, \underline{s}_{2}, \overline{s}_{1}, \overline{s}_{2}) := \esssup_{(l_1,l_2,m_1,m_2)\in \mathbb{A} (x,y_1,y_2,\underline{s}_{1}, \underline{s}_{2}, \overline{s}_{1}, \overline{s}_{2})} \mathbb{E} \left[w_T(x+ \right.  \nonumber \\
&& m_1 \underline{s}_{1} - l_1 \overline{s}_{1}+m_2 \underline{s}_{2}- l_2 \overline{s}_{2}, y_1-m_1+l_1,y_2-m_2-l_2,\underline{S}_{1}(T), \nonumber \\
&& \left. \underline{S}_{2}(T), \overline{S}_{1}(T), \overline{S}_{2}(T))|{\cal F}_{T-1}\right]
\end{eqnarray}
and for $k=1,2,\ldots,T$
\begin{eqnarray}\label{bel3}
&& w_{T-k}(x,y_1,y_2,\underline{s}_{1}, \underline{s}_{2}, \overline{s}_{1}, \overline{s}_{2}) := \esssup_{(l_1,l_2,m_1,m_2)\in \mathbb{A} (x,y_1,y_2,\underline{s}_{1}, \underline{s}_{2}, \overline{s}_{1}, \overline{s}_{2})} \mathbb{E} \left[w_{T-k+1}(x+ \right.  \nonumber \\
&& m_1 \underline{s}_{1} - l_1 \overline{s}_{1}+m_2 \underline{s}_{2} - l_2 \overline{s}_{2}, y_1-m_1+l_1,y_2-m_2+l_2,\underline{S}_{1}(T-k+1),\nonumber \\
&&\left. \underline{S}_{2}(T-k+1),
\overline{S}_{1}(T-k+1),\overline{S}_{2}(T-k+1))|{\cal F}_{T-k}\right].
\end{eqnarray}
In what follows to guarantee existence and continuity of random functions $w_t$ we shall assume that
\begin{assumption}\label{ass2}
\begin{equation}\label{Ass2}
 \mathbb{E}\left\{\overline{S}_{i}(t)\right\}<\infty,  \quad \mbox{and $\underline{S}_{i}(t)<\overline{S}_{i}(t)$ for $i=1,2$,  $t=1,2\ldots,T$}.
 \end{equation}
\end{assumption}

 Define for $k=1,2,\ldots,T$ and $(x,y_1,y_2)\in\mathbb{R}_{+} \times \mathbb{R}_{+}^{2}$
\begin{eqnarray}\label{bel4}
&&\bar{w}_{T-k}(x,y_1,y_2):= \mathbb{E}\left[w_{T-k+1}(x,y_1,y_2,\underline{S}_{1}(T-k+1), \underline{S}_{2}(T-k+1), \right. \nonumber \\
&& \left. \overline{S}_{1}(T-k+1),\overline{S}_{2}(T-k+1))|{\cal F}_{T-k}\right].
\end{eqnarray}
\begin{proposition}\label{submain}
Under Assumptions \ref{ass1}, \ref{ass2} functions $w_t$ and $\bar{w}_t$ are well defined for $t=0,1,\ldots, T$.
Moreover $\bar{w}_t$ is strictly concave. Furthermore there exists a unique optimal investment strategy for $w_t$.
\end{proposition}
\begin{proof} Notice first that by Assumption \ref{ass2} and concavity of $U$ using also Lemma 7.3 of \cite{RS2}  functions $w_t$ and $\bar{w}_t$ are well defined. By \eqref{bel4} we have
\begin{equation}
\bar{w}_{T-1}(x,y_1,y_2):= \mathbb{E}\left[w_T(x,y_1,y_2,\underline{S}_{1}(T), \underline{S}_{2}(T),
\overline{S}_{1}(T), \overline{S}_{2}(T))|{\cal F}_{T-1}\right].
\end{equation}
By Lemma \ref{concaveutil} $\bar{w}_{T-1}$ is a strictly concave function. Then as in the proof of Lemma \ref{concave} (see \eqref{con}) we have that the mapping
\begin{equation}\label{bel5}
(l_1,l_2,m_1,m_2)\mapsto \bar{w}_{T-1}(x+m_1 \underline{s}_1 - l_1 \overline{s}_1 + m_2 \underline{s}_2 - l_2 \overline{s}_2,y_1-m_1+l_1,y_2-m_2+l_2)
\end{equation}
is strictly concave within the class of controls $(l_1,l_2,m_1,m_2)$ such that $l_1m_1+l_2m_2=0$. Consequently as in Proposition \ref{Prop2.3} for each $(\underline{s}_1,\underline{s}_2,\overline{s}_1,\overline{s}_2)$ there is a unique optimal control $(\hat{l}_1,\hat{l}_2,\hat{m}_1,\hat{m}_2)$ for which $\esssup$ in \eqref{bel2} is attained and it is a continuous functions of $(x,y_1,y_2,\underline{s}_1,\underline{s}_2,\overline{s}_1,\overline{s}_2)$.
Assume that for $(x^1,y_1^1,y_1^1)$ the optimal strategy is $(\hat{l}_1^1,\hat{l}_2^1,\hat{m}_1^1,\hat{m}_2^1)$ while for  $(x^2,y_1^2,y_2^2)$ the optimal strategy is $(\hat{l}_1^2,\hat{l}_2^2,\hat{m}_1^2,\hat{m}_2^2)$. Then as in \eqref{con} we have
\begin{eqnarray}\label{bel6}
&&w_{T-1}(\lambda(x^1,y_1^1,y_2^1)+(1-\lambda)(x^2,y_1^2,y_2^2),\underline{s}_{1}, \underline{s}_{2}, \overline{s}_{1}, \overline{s}_{2})\geq  \nonumber \\
&&\bar{w}_{T-1}(\lambda (x^1+\hat{m}_1^1 \underline{s}_1 - \hat{l}_1^1 \overline{s}_1 + \hat{m}_2^1 \underline{s}_2 - \hat{l}_2^1 \overline{s}_2,y_1^1-\hat{m}_1^1+\hat{l}_1^1,y_2^1-\hat{m}_2^1+\hat{l}_2^1)+ \nonumber \\
&&(1-\lambda)
(x^2+\hat{m}_1^2 \underline{s}_1 - \hat{l}_1^2 \overline{s}_1 + \hat{m}_2^2 \underline{s}_2 - \hat{l}_2^2 \overline{s}_2,y_1^2-\hat{m}_1^2+\hat{l}_1^2,y_2^1-\hat{m}_2^2+\hat{l}_2^2))\geq \nonumber \\
&&\lambda w_{T-1}(x^1,y_1^1,y_2^1,\underline{s}_{1}, \underline{s}_{2}, \overline{s}_{1}, \overline{s}_{2})+
(1-\lambda) w_{T-1}(x^2,y_1^2,y_2^2,\underline{s}_{1}, \underline{s}_{2}, \overline{s}_{1}, \overline{s}_{2})
\end{eqnarray}
with strict inequality in the third line when
\begin{eqnarray}\label{bel7}
&&z_{T-1}(x^1,y_1^1,y_2^1,\underline{s}_1,\underline{s}_2,\overline{s}_1,\overline{s}_2):= \nonumber \\
&&(x^1+\hat{m}_1^1 \underline{s}_1 - \hat{l}_1^1 \overline{s}_1 + \hat{m}_2^1 \underline{s}_2 - \hat{l}_2^1 \overline{s}_2,y_1^1-\hat{m}_1^1+\hat{l}_1^1,y_2^1-\hat{m}_2^1+\hat{l}_2^1)\neq \nonumber \\
&&z_{T-1} (x^2,y_1^2,y_2^2,\underline{s}_1,\underline{s}_2,\overline{s}_1,\overline{s}_2)= \nonumber \\
&&(x^2+\hat{m}_1^2 \underline{s}_1 - \hat{l}_1^2 \overline{s}_1 + \hat{m}_2^2 \underline{s}_2 - \hat{l}_2^2 \overline{s}_2,y_1^2-\hat{m}_1^2+\hat{l}_1^2,y_2^2-\hat{m}_2^2+\hat{l}_2^2).
\end{eqnarray}
Furthermore when $\hat{l}_1^1 \hat{m}_1^2>0$ by optimality of $w_{T-1}$ we have
\begin{eqnarray}
&&w_{T-1}(x^1,y_1^1,y_2^1,\underline{s}_{1}, \underline{s}_{2}, \overline{s}_{1}, \overline{s}_{2})= \nonumber \\ &&\bar{w}_{T-1}(x^1 - \hat{l}_1^1 \overline{s}_1 + \hat{m}_2^1 \underline{s}_2 - \hat{l}_2^1 \overline{s}_2,y_1^1+\hat{l}_1^1,y_2^1-\hat{m}_2^1+\hat{l}_2^1)=\nonumber \\
&&w_{T-1}(\bar{x}^1-\hat{l}_1^1 \overline{s}_1,y_1^1+\hat{l}_1^1,\bar{y}_2^1,\underline{s}_{1}, \underline{s}_{2}, \overline{s}_{1}, \overline{s}_{2})=w_{T-1}(z^1)
\end{eqnarray}
with $\bar{x}^1=x^1+ \hat{m}_2^1 \underline{s}_2 - \hat{l}_2^1 \overline{s}_2$ and $\bar{y}_2^1=y_2^1-\hat{m}_2^1+\hat{l}_2^1$ and where we shall use a shorter notation $z^i= z_{T-1}(x^i,y_1^i,y_2^i,\underline{s}_1,\underline{s}_2,\overline{s}_1,\overline{s}_2)$ for $i=1,2$.
Similarly
\begin{eqnarray}
&&w_{T-1}(x^2,y_1^2,y_2^2,\underline{s}_{1}, \underline{s}_{2}, \overline{s}_{1}, \overline{s}_{2})= \nonumber \\ &&\bar{w}_{T-1}(x^2 + \hat{m}_1^2 \underline{s}_1 + \hat{m}_2^2 \underline{s}_2 - \hat{l}_2^2
 \overline{s}_2,y_1^2-\hat{m}_1^2,y_2^2-\hat{m}_2^2+\hat{l}_2^2)=\nonumber \\
&&w_{T-1}(\bar{x}^2+\hat{m}_1^2 \underline{s}_1,y_1^2-\hat{m}_1^2,\bar{y}_2^2,\underline{s}_{1}, \underline{s}_{2}, \overline{s}_{1}, \overline{s}_{2})=w_{T-1}(z^2)
\end{eqnarray}
with $\bar{x}^2=x^2 + \hat{m}_2^2 \underline{s}_2 - \hat{l}_2^2 \overline{s}_2$ and $\bar{y}_2^2=y_2^2-\hat{m}_2^2+\hat{l}_2^2$.
Then for $\lambda={\hat{m}_1^2 \over \hat{l}_1^1+\hat{m}_1^2}$ we have that $\lambda\hat{l}_1^1+(1-\lambda)\hat{m}_1^2=0$, and  when $z^1=z^2$  we obtain
\begin{eqnarray}\label{bel6p}
&&w_{T-1}(\lambda(x^1,y_1^1,y_2^1)+(1-\lambda)(x^2,y_1^2,y_2^2),\underline{s}_{1}, \underline{s}_{2}, \overline{s}_{1}, \overline{s}_{2})\geq \nonumber \\
&&w_{T-1}(\lambda(\bar{x}^1,y_1^1,\bar{y}_2^1)+(1-\lambda)(\bar{x}^2,y_1^2,\bar{y}_2^2),\underline{s}_{1}, \underline{s}_{2}, \overline{s}_{1}, \overline{s}_{2})\geq  \nonumber \\
&&\bar{w}_{T-1}(\lambda(\bar{x}^1,y_1^1,\bar{y}_2^1)+(1-\lambda)(\bar{x}^2,y_1^2,\bar{y}_2^2))>\bar{w}_{T-1}(\lambda \bar{x}^1 + (1-\lambda)\bar{x}^2 - \lambda\hat{l}_1^1 (\overline{s}_1-\underline{s}_1), \nonumber \\
&&\lambda y_1^1+ (1-\lambda) y_1^2, \lambda \bar{y}_2^1+ (1-\lambda)\bar{y}_2^2))=
\bar{w}_{T-1}(\lambda (\bar{x}^1 -\hat{l}_1^1\overline{s}_1) +  \nonumber \\
&&(1-\lambda)(\bar{x}^2 + \hat{m}_1^2\underline{s}_1), \lambda(y_1^1+\hat{l}_1^1)+ (1-\lambda)(y_1^2-\hat{m}_1^2), \lambda \bar{y}_2^1+ (1-\lambda)\bar{y}_2^2)= \nonumber\\
&&w_{T-1}(z^1)=w_{T-1}(z^2)= \lambda w_{T-1}(x^1,y_1^1,y_2^1,\underline{s}_{1}, \underline{s}_{2}, \overline{s}_{1}, \overline{s}_{2})+ \nonumber \\
&&(1-\lambda) w_{T-1}(x^2,y_1^2,y_2^2,\underline{s}_{1}, \underline{s}_{2}, \overline{s}_{1}, \overline{s}_{2})
\end{eqnarray}
Consequently we have strict inequality in \eqref{bel6}. Since $w_{T-1}$ is concave (as a function of the first three coordinates) and we have strict inequality in \eqref{bel6p} for a particular $\lambda\in [0,1]$  we have strict inequality in \eqref{bel6p} for any $\lambda\in (0,1)$. We have the same strict inequality in \eqref{bel6p} in the cases when $\hat{m}_1^1 \hat{l}_1^2>0$, $\hat{l}_2^1 \hat{m}_2^2>0$ and $\hat{m}_2^1 \hat{l}_2^2>0$ assuming  that $z^1=z^2$ (using similar consideration).

Let
\begin{eqnarray}\label{bel8}
&&\bar{w}_{T-2}(x,y_1,y_2):= \mathbb{E}\left[w_{T-1}(x,y_1,y_2,\underline{S}_{1}(T-1), \underline{S}_{2}(T-1), \right. \nonumber \\
&& \left. \overline{S}_{1}(T-1), \overline{S}_{2}(T-1))|{\cal F}_{T-2}\right].
\end{eqnarray}
Then by \eqref{bel6} we have
\begin{eqnarray}\label{bel85}
&&\bar{w}_{T-2}(\lambda(x^1,y_1^1,y_2^1)+(1-\lambda)(x^2,y_1^2,y_2^2))=
\mathbb{E}\left[w_{T-1}(\lambda(x^1,y_1^1,y_2^1)+\right.
\nonumber \\
&&\left. (1-\lambda)(x^2,y_1^2,y_2^2), \underline{S}_{1}(T-1), \underline{S}_{2}(T-1), \overline{S}_{1}(T-1), \overline{S}_{2}(T-1))|{\cal F}_{T-2}\right]\geq \nonumber \\
&& \mathbb{E}\left[\bar{w}_{T-1}(\lambda (x^1+\hat{m}_1^1 \underline{S}_1(T-1) - \hat{l}_1^1 \overline{S}_1(T-1) +\hat{m}_2^1 \underline{S}_2(T-1) -\right. \nonumber \\
&&  \hat{l}_2^1 \overline{S}_2(T-1),y_1^1-\hat{m}_1^1+\hat{l}_1^1,y_2^1-\hat{m}_2^1+\hat{l}_2^1))+(1-\lambda)
(x^2+\hat{m}_1^2 \underline{S}_1(T-1) - \nonumber \\
&& \left. \hat{l}_1^2 \overline{S}_1(T-1)+ \hat{m}_2^2 \underline{S}_2(T-1) - \hat{l}_2^2 \overline{S}_2(T-1),y_1^2-\hat{m}_1^2+\hat{l}_1^2, y_2^2-\hat{m}_2^2+\hat{l}_2^2))|{\cal F}_{T-2}\right]\geq\nonumber \\
&& \lambda \mathbb{E}\left[\bar{w}_{T-1}(x^1+\hat{m}_1^1 \underline{S}_1(T-1) - \hat{l}_1^1 \overline{S}_1(T-1)+ \hat{m}_2^1 \underline{S}_2(T-1) - \hat{l}_2^1 \overline{S}_2(T-1),\right.  \nonumber \\
&&\left.y_1^1-\hat{m}_1^1+\hat{l}_1^1,y_2^1-\hat{m}_2^1+\hat{l}_2^1))|{\cal F}_{T-2}\right]+ (1-\lambda) \mathbb{E} \left[\bar{w}_{T-1}(x^2+\hat{m}_1^2 \underline{S}_1(T-1) \right.\nonumber \\
&& \left. - \hat{l}_1^2 \overline{S}_1(T-1) + \hat{m}_2^2 \underline{S}_2(T-1) - \hat{l}_2^2 \overline{S}_2(T-1),y_1^2-\hat{m}_1^2+\hat{l}_1^2,y_2^2-\hat{m}_2^2+\hat{l}_2^2))|{\cal F}_{T-2}\right]\nonumber\\
&& = \lambda \mathbb{E}\left[w_{T-1}(x^1,y_1^1,y_2^1,\underline{S}_{1}(T-1), \underline{S}_{2}(T-1), \overline{S}_{1}(T-1), \overline{S}_{2}(T-1))|{\cal F}_{T-2}\right]+ \nonumber \\
&&(1-\lambda) \mathbb{E}\left[{w}_{T-1}(x^2,y_1^2,y_2^2,\underline{S}_{1}(T-1), \underline{S}_{2}(T-1), \overline{S}_{1}(T-1), \overline{S}_{2}(T-1))|{\cal F}_{T-2}\right].
\end{eqnarray}
We have a strict inequality on the ${\cal F}_{T-1}$ measurable set $Z_{T-1}\cap D_{T-1}$ where
\begin{equation}
D_{T-1}=\left\{\hat{l}_1^1 \hat{m}_1^2>0\right\}\cup\left\{\hat{m}_1^1 \hat{l}_1^2>0\right\}\cup \left\{\hat{l}_2^1 \hat{m}_2^2>0\right\} \cup \left\{\hat{m}_2^1 \hat{l}_2^2>0\right\}
\end{equation}
\begin{eqnarray}\label{bel9}
&& Z_{T-1}=\left\{z_{T-1}(x^1,y_1^1,y_2^1,\underline{S}_1(T-1),\underline{S}_2(T-1),\overline{S}_1(T-1),
\overline{S}_2(T-1))= \right. \nonumber \\
&& \left. z_{T-1} (x^2,y_1^2,y_2^2,\underline{S}_1(T-1),\underline{S}_2(T-1),\overline{S}_1(T-1),\overline{S}_2(T-1))\right\}.
\end{eqnarray}
On the set $Z_{T-1}$ we have
\begin{equation}\label{bel10}
x^1-x^2 + (\hat{m}_1^1-\hat{m}_1^2)\underline{S}_1(T-1) - (\hat{l}_1^1-\hat{l}_1^2) \overline{S}_1(T-1) + (\hat{m}_2^1-\hat{m}_2^2)\underline{S}_2(T-1) - (\hat{l}_2^1-\hat{l}_2^2) \overline{S}_2(T-1)=0
\end{equation}

We are going to show that $\mathbb{P}\left[Z_{T-1}\setminus D_{T-1}|{\cal F}_{T-2}\right]=0$.
Now notice that $\hat{l}_1^1\hat{m}_1^1=0$, $\hat{l}_2^1\hat{m}_2^1=0$, $\hat{l}_1^2\hat{m}_1^2=0$ and $\hat{l}_2^2 \hat{m}_2^2=0$. Assume that when $\hat{l}_i^1>0$ then $\hat{m}_i^2=0$ and when $\hat{m}_i^1>0$ then $\hat{l}_i^1=0$, for $i=1,2$. We have to consider 4 cases:

(i) $\hat{l}_1^1>0$, $\hat{m}_1^2=0$, $\hat{l}_2^1>0$ and $\hat{m}_2^2=0$,

(ii) $\hat{m}_1^1>0$, $\hat{l}_1^2=0$,  $\hat{l}_2^1>0$ and $\hat{m}_2^2=0$,

(iii) $\hat{l}_1^1>0$, $\hat{m}_1^2=0$, $\hat{m}_2^1>0$ and $\hat{l}_2^2=0$,

(iv) $\hat{m}_1^1>0$, $\hat{l}_1^2=0$, $\hat{m}_2^1>0$ and $\hat{l}_2^2=0$'

It is clear that it suffices to consider only the cases (i) and (ii) since (iii) and (iv) are similar to (ii) and (i).
\noindent
In the case  (i) on the set $Z_{T-1}$ we have
\begin{eqnarray}\label{bel10p}
&&x^1-x^2  - (\hat{l}_1^1-\hat{l}_1^2) \overline{S}_1(T-1) - (\hat{l}_2^1-\hat{l}_2^2) \overline{S}_2(T-1)= \nonumber \\
&& x^1-x^2 + (y_1^1-y_1^2) \overline{S}_1(T-1) + (y_2^1-y_2^2) \overline{S}_2(T-1) = 0
\end{eqnarray}
while in the case  (ii) on the set $Z_{T-1}$ we have
\begin{eqnarray}\label{bel10pp}
&& x^1-x^2 + (\hat{m}_1^1-\hat{m}_1^2)\underline{S}_1(T-1)  - (\hat{l}_2^1-\hat{l}_2^2) \overline{S}_2(T-1)= \nonumber \\
&&  x^1-x^2 + (y_1^1-y_1^2)\underline{S}_1(T-1) + (y_1^1-y_1^2)\overline{S}_2(T-1)=0
\end{eqnarray}
Consequently by Assumption \ref{ass1} we have that $\mathbb{P}\left[Z_{T-1}\setminus D_{T-1}|{\cal F}_{T-2}\right]=0$ and $\bar{w}_{T-2}$ is a strictly concave function.
We now continue the proof by induction. Given strict concavity of function $\bar{w}_{T-k+1}$ following \eqref{bel85}-\eqref{bel10pp} we obtain strict concavity of $\bar{w}_{T-k}$. Finally by Proposition \ref{Prop2.3} we have existence of unique optimal investment strategies.
\end{proof}
We are now in position to formulate main result of the paper
\begin{theorem} Under the assumptions of Proposition \ref{submain} we have that
\begin{equation}\label{mainr}
w_0(x,y_1,y_2,\underline{S}_{1}(0), \underline{S}_{2}(0), \overline{S}_{1}(0), \overline{S}_{2}(0))=\sup E\left\{U(W_T)\right\}
\end{equation}
where the suppremum is taken over all admissible strategies with $\underline{S}_{1}(t), \underline{S}_{2}(t), \overline{S}_{1}(t), \overline{S}_{2}(t)$ for $t=0,1,\ldots,T$ as bid and ask prices.
Furthermore $w_0$ is also a supremum for the shadow prices defined in Theorem \ref{sshadow} in section 8 for each time $t\in \left\{0,1,\ldots,T\right\}$.
\end{theorem}
\begin{proof}
To prove \eqref{mainr} we use standard arguments of the paper \cite{RS2} using properties of regular conditional probability shown in Lemma 7.3 of \cite{RS2}. The second part of Theorem follows directly from the static result shown in Theorem \ref{sshadow} in section 8.
\end{proof}

\section{Generalization to multiasset case}
We consider now the case of several assets as was mentioned in the introduction. We point out only major steps since the method is an easy extension of the case $d=2$ studied earlier in details. We consider first static model consisting of bid $\underline{S}=(\underline{s}_i)_{i=1,\ldots,d}$ and ask prices $\overline{S}=(\overline{s}_i)_{i=1,\ldots,d}$ and the set of prices 

\begin{equation}
    \mathbb{D} := \big\{ (\underline{s}_{i}), (\overline{s}_{i}) \in \mathbb{R}_{+}^{d} \times \mathbb{R}_{+}^{d} : \quad 0 < \underline{s}_{i} \leq \overline{s}_{i} \quad \mbox{for} \quad i=1,\ldots,d\big\}.
\end{equation}
Financial position will be denoted by $(x, (y_{i})_{i=1,\ldots,d})$ and will be nonnegative.
For every $(x,(y_{i}),(\underline{s}_{i}), (\overline{s}_{i})) \in \mathbb{R}_{+} \times \mathbb{R}_{+}^{d} \times \mathbb{D}$ we have the following set of admissible investment strategies
\begin{eqnarray}
&&\mathbb{A}((x,(y_{i}),(\underline{s}_{i}), (\overline{s}_{i})) := \Big\{ ((l_{i}), (m_{i})) \in \mathbb{R}_{+}^{d} \times \mathbb{R}_{+}^{d} :\nonumber \\
&& \quad 0 \leq x + \sum_{i=1}^d(\underline{s}_{i} m_{i} - \overline{s}_{i}l_{i}),  0 \leq y_{i} - m_{i} + l_{i}, m_{i} \leq y_{i} \quad \mbox{for} \quad i=1,\ldots,d \Big\} .
\end{eqnarray}
Let
\begin{equation}
    \begin{split}
        \mathbb{P} := \Big\{ (x, (y_{i}),(\underline{s}_{i}), (\overline{s}_{i}), (l_{i}), & (m_{i})) \in \mathbb{R}_{+} \times \mathbb{R}_{+}^{d} \times \mathbb{D} \times \mathbb{R}_{+}^{d} \times \mathbb{R}_{+}^{d}:  \\
        & ((l_{i}), (m_{i})) \in \mathbb{A}(x, (y_{i}), (\underline{s}_{i}), (\overline{s}_{i})) \Big\} .
    \end{split}
\end{equation}

Let $g : \mathbb{R}_{+} \times \mathbb{R}_{+}^{d} \longrightarrow \mathbb{R} \cup \{ - \infty \}$ be a strictly concave function which is strictly increasing with respect to each variable. 

Define the function $\varphi : \mathbb{P} \longrightarrow \mathbb{R} \cup \{ - \infty \}$ by the formula
\begin{equation}\label{vphim}
\varphi (x, (y_{i}), (\underline{s}_{i}), (\overline{s}_{i}), (l_{i}), (m_{i})) = g \Big( x + \sum_{i=1}^d (\underline{s}_{i} m_{i} - \overline{s}_{i} l_{i}), (y_{i} - m_{i}+l_{i})  \Big)
\end{equation}
and
$\Phi : \mathbb{R}_{+} \times \mathbb{R}_{+}^{d} \times \mathbb{D} \longrightarrow \mathbb{R} \cup \{ - \infty \}$ for every $(x, (y_{i}), (\underline{s}_{i}), (\overline{s}_{i}), (l_{i}), (m_{i})) \in \mathbb{R}_{+} \times \mathbb{R}_{+}^{d} \times \mathbb{D}$
\begin{eqnarray}\label{Phim}
&&\Phi (x, (y_{i}), (\underline{s}_{i}), (\overline{s}_{i}), (l_{i}), (m_{i})) := \nonumber \\
&&\sup_{((l_{i}), (m_{i})) \in \mathbb{A}(x, (y_{i}), (\underline{s}_{i}), (\overline{s}_{i}), (l_{i}), m_{i}))} \varphi (x, (y_{i}), (\underline{s}_{i}), (\overline{s}_{i}), (l_{i}), (m_{i}))   .
\end{eqnarray}

For every $ (x, (y_{i}), (\underline{s}_{i}), (\overline{s}_{i}))\in \mathbb{R}_{+} \times \mathbb{R}_{+}^{d} \times \mathbb{D}$ define the set
    \begin{eqnarray}
        && \overleftrightarrow{\mathbb{A}}(x, (y_{i}), (\underline{s}_{i}), (\overline{s}_{i}) :=  \nonumber \\
        && \Big\{ (l_{1}, l_{2}, m_{1}, m_{2}) \in \mathbb{A} (x, (y_{i}), (\underline{s}_{i}), (\overline{s}_{i})),\quad \sum_{i=1}^d l_{i} m_{i} = 0 \Big\}.
    \end{eqnarray}
By analogy to 
Proposition \ref{continuity of widtilde} 
and Corollary \ref{Corollary on the continuity of the function widehat}
as well as Theorem \ref{select} we have 

\begin{proposition}\label{continuitym} 
There exists a function  $u : \mathbb{R}_{+} \times \mathbb{R}_{+}^{d} \times \mathbb{D} \rightarrow \mathbb{R}_{+} \times \mathbb{R}_{+}^{d}$ be such that for every $(x, (y_{i}), (\underline{s}_{i}), (\overline{s}_{i})) \in \mathbb{R}_{+} \times \mathbb{R}_{+}^{d} \times \mathbb{D}$ we have that
    \begin{equation*}
        u(x, (y_{i}), (\underline{s}_{i}), (\overline{s}_{i})) \in \overleftrightarrow{\mathbb{A}}(x, (y_{i}), (\underline{s}_{i}), (\overline{s}_{i}))
    \end{equation*}
    and
    \begin{equation*}
        \Phi (x, (y_{i}), (\underline{s}_{i}), (\overline{s}_{i})) = \varphi \big(x, (y_{i}), (\underline{s}_{i}), (\overline{s}_{i}), u(x, (y_{i}), (\underline{s}_{i}), (\overline{s}_{i})) \big).
    \end{equation*}
Furthermore $u$ is continuous and uniquely defined. Additionally, for fixed $j\in \left\{1,\ldots,d\right\}$ consider  $\widehat{u} : [\underline{s}_{j}, \overline{s}_{j}] \longrightarrow \mathbb{R}_{+}^{d} \times \mathbb{R}_{+}^{d}$ such that 
         $  \widehat{u}(s) \in \overleftrightarrow{\mathbb{A}}(x, (y_{i}), (\check{\underline{s}}_{i}(s)), (\check{\overline{s}}_{i}(s)))$ for every $s \in [\underline{s}_{j}, \overline{s}_{j}]$ and
    \begin{equation*}
        \Phi (x, (y_{i}), y_{2},(\check{\underline{s}}_{i}(s)),(\check{\overline{s}}_{i}(s)) ) = \varphi \big(x, (y_{i}), y_{2},(\check{\underline{s}}_{i}(s)),(\check{\overline{s}}_{i}(s)), \widehat{u}(s) \big), 
         \end{equation*}
    where $\check{\underline{s}}_{i}(s)=s$ for $i=j$ and $\underline{s}_{i}$ otherwise and similarly 
    $\check{\overline{s}}_{i}(s)=s$ for $i=j$ and $\overline{s}_{i}$ otherwise. 
             Then $\widehat{u}$ is a  continuous function.
\end{proposition}
For every $(i_{1}, i_{2},\ldots, i_{d}) \in \{ -1, 0, 1 \}^{d}$ and for every $x, (y_{i}), (\underline{s}_{i}), (\overline{s}_{i}) \in \mathbb{R}_{+} \times \mathbb{R}_{+}^{d} \times \mathbb{D}$  define
\begin{eqnarray}
&& \mathbb{A}^{(i_{1}, i_{2},\ldots,i_{d})} (x, (y_{i}), (\underline{s}_{i}), (\overline{s}_{i}))
        := \Big\{ (l_{i}, m_{i}) \in \mathbb{A}(x, (y_{i}), (\underline{s}_{i}), (\overline{s}_{i})) : \nonumber \\
&&\quad i_{j} = \sgn (l_{j} - m_{j}) \quad \mbox{for} \quad j=1,\ldots,d, \quad \sum_{i=1}^d l_{i} m_{i} = 0 \Big\}.
\end{eqnarray}
We also define trading sets for every $(i_{1}, i_{2},\ldots, i_{d}) \in \{ -1, 0, 1 \}^{d}$ 
\begin{eqnarray}
    &&\mathbf{t}^{(i_{1}, i_{2},\ldots, i_{d})}   ((\underline{s}_{i}), (\overline{s}_{i})) := \Big\{ (x, y_{1}, y_{2}) \in \mathbb{R}_{+} \times \mathbb{R}_{+}^{2} : \nonumber \\
    &&\exists_{((l_{i}), (m_{i})) \in \mathbb{A}^{(i_{1}, i_{2},\ldots, i_{d})} (x, (y_{i}), (\underline{s}_{i}), (\overline{s}_{i}))} \Phi (x, (y_{i}), (\underline{s}_{i}), (\overline{s}_{i})) =  \nonumber \\
    && \varphi (x, (y_{i}), (\underline{s}_{i}), (\overline{s}_{i})), (l_{i}), (m_{i})) \Big\}.
\end{eqnarray}
and in analogy to Lemma \ref{disjoint} we have
\begin{lemma}\label{disjointm}
 $\Big\{ \mathbf{t}^{(i_{1}, i_{2}, \ldots,i_{d})} ( (\underline{s}_{i}), (\overline{s}_{i})) \Big\}_{(i_{1}, i_{2},\ldots,i_{d}) \in \{ - 1, 0, 1 \}}$ is a family of disjoint Borel sets which cover $\mathbb{R}_{+} \times \mathbb{R}_{+}^{d}$.    
\end{lemma}
Finally to complete multidimensional static case we need (see  Theorem \ref{Thm7.1})
\begin{theorem}\label{Thm7.1m}
 For   $((\underline{s}_{i}), (\overline{s}_{i})) \in \mathbb{D}$ and
 $(i_2,i_3\ldots,i_d) \in \{ - 1, 0, 1 \}^{d-1}$ we have
 \begin{equation}\label{nodecom}
  \mathbf{t}^{(0, i_2,i_3,\ldots,i_d)}((\underline{s}_{i}), (\overline{s}_{i})) 
 = \bigcup_{s\in [\underline{s}_{1},\overline{s}_{1}]} \mathbf{t}^{(0, i_2,i_3,\dots,i_d)}({s}, \underline{s}_{2},\ldots, \underline{s}_d, {s}, \overline{s}_{2},\ldots, \overline{s}_{d}), 
\end{equation} 
and similarly with $i_j=0$ and buying and selling $j$-th asset price $s$.
\end{theorem}
Let $(\underline{s}_{i},\overline{s}_{i}) \in \mathbb{D}$.
From Lemma \ref{disjoint} we know that for every $(x, (y_{i})) \in \mathbb{R}_{+} \times \mathbb{R}_{+}^{d}$ there exists a unique sequence $(i_{1}, i_{2},\ldots,i_d) \in \{ - 1, 0, 1 \}^{d}$ such that $(x, (y_{i})) \in \mathbf{t}^{(i_{1}, i_{2},\ldots,i_d)}((\underline{s}_{i}), (\overline{s}_{i}))$.
Consequently for  $(x, (y_{i})) \in \mathbb{R}_{+} \times \mathbb{R}_{+}^{d}$ we can define
$\mathrm{I}_{j}(x, (y_{i})) \in \{ - 1, 0, 1 \}$
such that for $j=1,2,\ldots,d$
\begin{equation}
    (x, (y_{i})) \in \mathbf{t}^{\big( \mathrm{I}_{1}(x, (y_{i})), \mathrm{I}_{2}(x, (y_{i})), \dots, \mathrm{I}_{d}(x, (y_{i})) \big)}((\underline{s}_{i}), (\overline{s}_{i})).
\end{equation}

From Theorem \ref{Thm7.1m}  for every $(i_2,i_3\ldots,i_d) \in \{ - 1, 0, 1 \}^{d-1}$ and for every $(x, (y_{i}), (y_{i})) \in \mathbf{t}^{(0, i_2,\ldots,i_d)}((\underline{s}_{i}), (\overline{s}_{i}))$ there exists an $\widetilde{s} \in [\underline{s}_{1}, \overline{s}_{1}]$ such that 
\begin{equation*}
(x, (y_{i})) \in \mathbf{t}^{(0,i_2,\ldots,i_d)}(\widetilde{s}, \underline{s}_{2},\ldots, \underline{s}_d, \widetilde{s}, \overline{s}_{2},\ldots,\overline{s}_d).
\end{equation*}
Therefore, for every $(x, (y_{i})) \in \mathbb{R}_{+} \times \mathbb{R}_{+}^{2d}$ we can define the set

\begin{equation}\label{sha1m}
    {\cal  S}_{1}(x, (y_{i})) :=             \begin{cases}
                \{ \underline{s}_{1} \} \quad & \mbox{if $\mathrm{I}_{1}(x, (y_{i})) = - 1$} \\
                \Big\{ s \in [\underline{s}_{1}, \overline{s}_{1}] : & \mbox{when $(x, (y_{i})) \in \mathbf{t}^{\big( 0, \mathrm{I}_{2}(x, (y_{i})),\ldots,\mathrm{I}_{d}(x, (y_{i})) \big)} (s, \underline{s}_{2},\ldots,\underline{s}_{d},$} \\
                & \mbox{ $s, \overline{s}_{2},\ldots,\overline{s}_{d})$} \quad \mbox{and $\mathrm{I}_{1}(x, (y_{i})) = 0\Big\}$} \\
                \{ \overline{s}_{1} \} \quad & \mbox{if $\mathrm{I}_{1}(x, (y_{i})) = 1$}
            \end{cases}.
\end{equation}

By analogy we can also define sets ${\cal S}_j(x,(y_{i}))$ for $j=2,\ldots,d$. Then for  
$\widetilde{S}_j(x,(y_i)):=\inf{\cal S}_j(x,y(i))$   we have an analog of Theorem \ref{sshadow}

\begin{theorem}\label{sshadowm}
The sequence $\left(\widetilde{S}_{1}(x, (y_{i})),\ldots, \widetilde{S}_{d}(x, (y_{i}))\right)$ defined above form shadow price for the static problem i.e. with $\Phi$ defined in \eqref{Phim} we have
\begin{eqnarray}\label{shaddm}
&&\Phi (x, (y_{i}), (\underline{s}_{i}), (\overline{s}_{i})) ) := \nonumber \\
&&\sup_{((l_{i}), (m_{i})) \in \mathbb{A}(x, (y_{i}),\widetilde{S}_{1},\ldots, \widetilde{S}_{d}, \widetilde{S}_{1}, \dots, \widetilde{S}_{d})} \varphi (x, (y_{i}), \widetilde{S}_1,\ldots, \widetilde{S}_{d}, \widetilde{S}_{1}, \ldots \widetilde{S}_{d}, (l_{i}), (m_{i})).
\end{eqnarray}
where to simplify notations we write $\widetilde{S}_j$ for $\widetilde{S}_{j}(x, (y_{i}))$.
\end{theorem}

Consider regular conditional laws $\mathbb{P}_t^{i_1,\ldots,i_d}$ (which exist by Theorem 6.3 of \cite{Kal}) defined for $i_1,i_2,\ldots,i_d\in\left\{-1,1\right\}$ as
\begin{equation}\label{Ass1m}
\mathbb{P}_t^{i_1,i_2,\ldots,i_d}(\cdot)=\mathbb{P}\left((S_1^{i_1}(t+1),\ldots,S_d^{i_d}(t+1))\in \cdot|{\cal F}_t\right),
\end{equation}
with the notation $S_k^i(t+1)=\underline{S}_{k}(t+1)$ when $i=-1$ and $S_k^i(t+1)=\overline{S}_{k}(t+1)$ when $i=1$ with $k=1,2,\ldots,d$.
We shall assume that
\begin{assumption}\label{ass1m}
Conditional laws $\mathbb{P}_t^{i_1,i_2,\ldots,i_d}$  for $i_1,i_2,\ldots,i_d\in\left\{-1,1\right\}$ are non degenerate, which means that
\begin{equation}
\mathbb{P}\left(\sum_{k=1}^d a_k S_k^{i_k}(t+1) = a_0 |{\cal F}_t\right)=0
\end{equation}
for any ${\cal F}_t$ measurable random variables $a_0,a_1,a_2,\ldots,a_d$ and $t=0,1,\ldots,T-1$. 
\end{assumption}
Define now using regular conditional probabilities  the following system of random Bellman equations for $(x,(y_i),(\underline{s}_{i}),(\overline{s}_{i})) \in \mathbb{R}_{+} \times \mathbb{R}_{+}^{d} \times \mathbb{D}$

\begin{equation}\label{bel1m}
w_T(x,(y_i),(\underline{s}_{i}),(\overline{s}_{i})):=U(x+\sum_{i=1}^d y_i\underline{s}_{i}),
\end{equation}
\begin{eqnarray}\label{bel2m}
&& w_{T-1}(x,(y_i),(\underline{s}_{i}), (\overline{s}_{i})) := \esssup_{((l_i),(m_i))\in \mathbb{A} (x,(y_i),(\underline{s}_{i}), (\overline{s}_{i}))} \mathbb{E} \left[w_T(x+ \right.  \nonumber \\
&& \sum_{i=1}^d (m_i \underline{s}_{i} - l_i \overline{s}_{i}), y_1-m_1+l_1,y_2-m_2-l_2, \ldots,y_d-m_d-l_d, \underline{S}_{1}(T), \nonumber \\
&& \left. \underline{S}_{2}(T),\ldots,\underline{S}_{d}(T),  \overline{S}_{1}(T), \overline{S}_{2}(T),\ldots, \overline{S}_{d}(T))|{\cal F}_{T-1}\right]
\end{eqnarray}
and for $k=1,2,\ldots,T$
\begin{eqnarray}\label{bel3m}
&& w_{T-k}(x,(y_i),(\underline{s}_{i}), (\overline{s}_{i})) :=
\esssup_{((l_i),(m_i))\in \mathbb{A} (x,(y_i),(\underline{s}_{i}), (\overline{s}_{i})} \mathbb{E} \left[w_{T-k+1}(x+ \right.  \nonumber \\
&&\sum_{i=1}^d m_i \underline{s}_{i} - l_i \overline{s}_{i}, y_1-m_1+l_1,y_2-m_2+l_2, \ldots, y_d-m_d+l_d, \underline{S}_{1}(T-k+1),\nonumber \\
&&\left. \underline{S}_{2}(T-k+1), \ldots,  \underline{S}_{d}(T-k+1),
\overline{S}_{1}(T-k+1),\overline{S}_{2}(T-k+1),\ldots,\right. \nonumber \\
&& \left. \overline{S}_{d}(T-k+1) )|{\cal F}_{T-k}\right].
\end{eqnarray}
In what follows to guarantee existence and continuity of random functions $w_t$ we shall assume that
\begin{assumption}\label{ass2m}
\begin{equation}\label{Ass2m}
 \mathbb{E}\left\{\overline{S}_{i}(t)\right\}<\infty,  \quad \mbox{and $\underline{S}_{i}(t)<\overline{S}_{i}(t)$ for $i=1,2,\ldots,d$,  $t=1,2\ldots,T$}.
 \end{equation}
\end{assumption}

 Define for $k=1,2,\ldots,T$ and $(x,(y_i))\in\mathbb{R}_{+} \times \mathbb{R}_{+}^{d}$
\begin{eqnarray}\label{bel4m}
&&\bar{w}_{T-k}(x,(y_i)):= \mathbb{E}\left[w_{T-k+1}(x,(y_i),\underline{S}_{1}(T-k+1), \underline{S}_{2}(T-k+1), \ldots, \right. \nonumber \\
&& \left. \underline{S}_{d}(T-k+1),\overline{S}_{1}(T-k+1),\overline{S}_{2}(T-k+1),\ldots, \overline{S}_{d}(T-k+1) )|{\cal F}_{T-k}\right].
\end{eqnarray}
By analogy to Proposition \ref{submain} we have 
\begin{proposition}\label{submainm}
Under Assumptions \ref{ass1m}, \ref{ass2m} functions $w_t$ and $\bar{w}_t$ are well defined for $t=0,1,\ldots, T$.
Moreover $\bar{w}_t$ is strictly concave. Furthermore there exists a unique optimal investment strategy for $w_t$.
\end{proposition}
The main result is stated as follows
\begin{theorem} Under the assumptions of Proposition \ref{submainm} we have that
\begin{equation}\label{mainr}
w_0(x,(y_i),(\underline{S}_{i}(0)), (\overline{S}_{i}(0)))=\sup E\left\{U(W_T)\right\}
\end{equation}
where the supremum is taken over all admissible strategies with bid and ask prices repsectively $\underline{S}_{1}(t), \underline{S}_{2}(t), \ldots  \underline{S}_{d}(t), \overline{S}_{1}(t), \overline{S}_{2}(t), \ldots, \overline{S}_{d}(t)$ for $t=0,1,\ldots,T$.
Furthermore $w_0$ is also a supremum for the shadow prices defined in Theorem \ref{sshadowm}  for each time $t\in \left\{0,1,\ldots,T\right\}$.
\end{theorem}
\appendix
\section{Some useful results}

\begin{lemma}\label{Lem1.2}
    Let $0 < \underline{s} \leq \overline{s}$ and $(l, m) \in \mathbb{R}_{+} \times \mathbb{R}_{+}$. Then
    \begin{equation}\label{simple inequality on wasting}
        \underline{s} m - \overline{s} l \leq \underline{s} (m - l)^{+} - \overline{s} (l - m)^{+}
    \end{equation}
with strict inequality when $\underline{s} < \overline{s}$ and $l m\neq 0$.
\end{lemma}

We say that random vector $Z=(Z_1,Z_2,\ldots, Z_d)$  defined on a probability space $(\Omega, {\cal F},\mathbb{P})$ is nondegenerate if there are no deterministic values $a_0,a_1,a_2,\ldots,a_d\in \mathbb{R}$ such that $\sum_{i=1}^d a_i Z_i=a_0$, $\mathbb{P}$ a.e..

\begin{lemma}\label{concaveutil}
Assume $U: \mathbb{R}_+ \mapsto \mathbb{R}$ is a strictly concave function and $Z=(Z_1,Z_2,\ldots,Z_d)$ is a $d$ - dimensional nondegenerate  random vector taking positive values. Then the function $g:\mathbb{R_+}\times \mathbb{R}_+^{2} \to \mathbb{R}$ defined as follows
\begin{equation}
g(x,y_1,y_2):= E\left\{U(x+ \sum_{i=1}^d y_iZ_i\right\}
\end{equation}
is also strictly concave, assuming that it is well defined.
\end{lemma}

\begin{proof} We have to show that for $\alpha\in (0,1)$, $(x,y_1,y_2,\ldots,y_d)\in \mathbb{R_+}^{d+1}$ and $(x',y_1',y_2')\in \mathbb{R_+}^{d+1}$ such that $(x,y_1,y_2,\dots,y_d)\neq (x',y_1',y_2',\ldots,y_d')$ we have
\begin{eqnarray}
&&g\left(\alpha x + (1-\alpha) x',\alpha y_1 + (1-\alpha) y_1',\ldots,\alpha y_d + (1-\alpha) y_d'\right)> \nonumber \\
&& \alpha g(x,y_1,\ldots,y_d) + (1-\alpha) g(x',y_1',\ldots,y_d').
\end{eqnarray}
Since $U$ is strictly concave we have that
\begin{eqnarray}
&&E\left\{U(\alpha(x+\sum_{i=1}^d y_iZ_i)+(1-\alpha)(x'+\sum_{i=1}^d y_i'Z_i))\right\}\geq \nonumber \\
&&\alpha E\left\{U(x+\sum_{i=1}^dy_iZ_i\right\} + (1-\alpha) E\left\{U(x'+\sum_{i=1}^d y_i'Z_i)\right\}
\end{eqnarray}
with equality only when $x+\sum_{i=1}^d y_iZ_i(\omega)=x'+ \sum_{i=1}^d y_i'Z_i(\omega)$, $\mathbb{P}$ a.e. However this contradics nondegeneracy of the random vector $Z$.
\end{proof}

\begin{theorem}\label{select}
     Assume that $(X, d)$  and  $(Y, \rho)$ are metric spaces and  for every $x \in X$ the set $\mathcal{A}(x) \subseteq Y$ is nonempty and compact. Assume furthermore that the mapping $x \longmapsto \mathcal{A}(x)$ is continuous in Hausdorff metric and for continuous function $\beta : X \times Y \longrightarrow \overline{\mathbb{R}}$ and every $x \in X$ there exists a unique $\widehat{a}(x) \in \mathcal{A}(x)$ such that
    \begin{equation}\label{eq:sel}
        \sup_{a \in \mathcal{A}(x)} \beta(x, a) = \beta \big( x, \widehat{a}(x) \big) .
    \end{equation}
    Then the mapping $x \longmapsto \widehat{a}(x)$ is continuous.
\end{theorem}
\begin{proof} Assume that $X \ni x_n \to x_0$. Since the mapping $x \longmapsto \mathcal{A}(x)$ is continuous in Hausdorff metric the set
\begin{equation*}
        \closure \Bigg( \mathcal{A}(x_{0}) \cup \bigcup_{n = 0}^{\infty} \mathcal{A}(x_{n}) \Bigg)
\end{equation*}
is a compact subset of $Y$. Therefore there exists such a subsequence $(n_{k})$ such that the sequence $\big( \widehat{a}(x_{n_{k}}) \big)$ converges to an element $\widetilde{a} \in \mathcal{A}(x_{0})$. By continuity of of $\beta$ we have
\begin{equation*}
         \beta \big( x_{0}, \widehat{a}(x_{0}) \big) = \lim_{k \rightarrow \infty} \beta \big( x_{n_{k}}, \widehat{a}(x_{n_{k}}) \big) = \beta(x_{0}, \widetilde{a}) .
\end{equation*}
Since $\widehat{a}(x_{0})$ is unique we have that $\widehat{a}(x_{0})=\widetilde{a}$. Hence for every convergent subsequence $\big( \widehat{a}(x_{n_{k}}) \big)$  has the same limit $\widehat{a}(x_{0})$ which  completes the proof of continuity of $\widehat{a}(\cdot)$.
\end{proof}


\bibliographystyle{siamplain}


\end{document}